\input amstex
\documentstyle{amsppt}
\magnification=\magstep 1
\document
\catcode`\@=11
\def\nologo{\let\logo@=\relax}
\catcode`\@=\active
\nologo
\vsize6.7in



\chardef\oldatsign=\catcode`\@
\catcode`\@=11
\newif\ifdraftmode			
\global\draftmodefalse


%
\font@\twelverm=cmr12 
\font@\twelvei=cmmi12 \skewchar\twelvei='177 
\font@\twelvesy=cmsy10 scaled\magstep1 \skewchar\twelvesy='060 
\font@\twelveex=cmex10 scaled\magstep1 
\font@\twelvemsa=msam10 scaled\magstep1 
\font@\twelvemsb=msbm10 scaled\magstep1 
\font@\twelvebf=cmbx12 
\font@\twelvett=cmtt12 
\font@\twelvesl=cmsl12 
\font@\twelveit=cmti12 
\font@\twelvesmc=cmcsc10 scaled\magstep1 
%
%
\font@\ninerm=cmr9 
\font@\ninei=cmmi9 \skewchar\ninei='177 
\font@\ninesy=cmsy9 \skewchar\ninesy='60 
\font@\ninemsa=msam9
\font@\ninemsb=msbm9
\font@\ninebf=cmbx9
%
%
%
\font@\ttlrm=cmbx12 scaled \magstep2 
\font@\ttlsy=cmsy10 scaled \magstep3 
\font@\tensmc=cmcsc10 
%
%
\def\normaltype{
	\def\pointsize@{12}%
	\abovedisplayskip18\p@ plus5\p@ minus9\p@
	\belowdisplayskip18\p@ plus5\p@ minus9\p@
	\abovedisplayshortskip1\p@ plus3\p@
	\belowdisplayshortskip9\p@ plus3\p@ minus4\p@
	\textonlyfont@\rm\twelverm
	\textonlyfont@\it\twelveit
	\textonlyfont@\sl\twelvesl
	\textonlyfont@\bf\twelvebf
	\textonlyfont@\smc\twelvesmc
	\ifsyntax@
		\def\big##1{{\hbox{$\left##1\right.$}}}%
	\else
		\let\big\twelvebig@
 \textfont0=\twelverm \scriptfont0=\ninerm \scriptscriptfont0=\sevenrm
 \textfont1=\twelvei  \scriptfont1=\ninei  \scriptscriptfont1=\seveni
 \textfont2=\twelvesy \scriptfont2=\ninesy \scriptscriptfont2=\sevensy
 \textfont3=\twelveex \scriptfont3=\twelveex  \scriptscriptfont3=\twelveex
 \textfont\itfam=\twelveit \def\it{\fam\itfam\twelveit}%
 \textfont\slfam=\twelvesl \def\sl{\fam\slfam\twelvesl}%
 \textfont\bffam=\twelvebf \def\bf{\fam\bffam\twelvebf}%
 \scriptfont\bffam=\ninebf \scriptscriptfont\bffam=\sevenbf
 \textfont\ttfam=\twelvett \def\tt{\fam\ttfam\twelvett}%
 \textfont\msafam=\twelvemsa \scriptfont\msafam=\ninemsa
 \scriptscriptfont\msafam=\sevenmsa
 \textfont\msbfam=\twelvemsb \scriptfont\msbfam=\ninemsb
 \scriptscriptfont\msbfam=\sevenmsb
	\fi
 \normalbaselineskip=\twelvebaselineskip
 \setbox\strutbox=\hbox{\vrule height12\p@ depth6\p@
      width0\p@}%
 \normalbaselines\rm \ex@=.2326ex%
}
%
%
%
\def\smalltype{
	\def\pointsize@{10}%
	\abovedisplayskip12\p@ plus3\p@ minus9\p@
	\belowdisplayskip12\p@ plus3\p@ minus9\p@
	\abovedisplayshortskip\z@ plus3\p@
	\belowdisplayshortskip7\p@ plus3\p@ minus4\p@
	\textonlyfont@\rm\tenrm
	\textonlyfont@\it\tenit
	\textonlyfont@\sl\tensl
	\textonlyfont@\bf\tenbf
	\textonlyfont@\smc\tensmc
	\ifsyntax@
		\def\big##1{{\hbox{$\left##1\right.$}}}%
	\else
		\let\big\tenbig@
	\textfont0=\tenrm \scriptfont0=\sevenrm \scriptscriptfont0=\fiverm 
	\textfont1=\teni  \scriptfont1=\seveni  \scriptscriptfont1=\fivei
	\textfont2=\tensy \scriptfont2=\sevensy \scriptscriptfont2=\fivesy 
	\textfont3=\tenex \scriptfont3=\tenex \scriptscriptfont3=\tenex
	\textfont\itfam=\tenit \def\it{\fam\itfam\tenit}%
	\textfont\slfam=\tensl \def\sl{\fam\slfam\tensl}%
	\textfont\bffam=\tenbf \def\bf{\fam\bffam\tenbf}%
	\scriptfont\bffam=\sevenbf \scriptscriptfont\bffam=\fivebf
	\textfont\msafam=\tenmsa
	\scriptfont\msafam=\sevenmsa
	\scriptscriptfont\msafam=\fivemsa
	\textfont\msbfam=\tenmsb
	\scriptfont\msbfam=\sevenmsb
	\scriptscriptfont\msbfam=\fivemsb
		\textfont\ttfam=\tentt \def\tt{\fam\ttfam\tentt}%
	\fi
 \normalbaselineskip 14\p@
 \setbox\strutbox=\hbox{\vrule height10\p@ depth4\p@ width0\p@}%
 \normalbaselines\rm \ex@=.2326ex%
}

\def\titletype{
	\def\pointsize@{17}%
	\textonlyfont@\rm\ttlrm
	\ifsyntax@
		\def\big##1{{\hbox{$\left##1\right.$}}}%
	\else
		\let\big\twelvebig@
		\textfont0=\ttlrm \scriptfont0=\twelverm
		\scriptscriptfont0=\tenrm
		\textfont2=\ttlsy \scriptfont2=\twelvesy
		\scriptscriptfont2=\tensy
	\fi
	\normalbaselineskip 25\p@
	\setbox\strutbox=\hbox{\vrule height17\p@ depth8\p@ width0\p@}%
	\normalbaselines
	\rm
	\ex@=.2326ex%
}

\def\tenbig@#1{
	{%
		\hbox{%
			$%
			\left
			#1%
			\vbox to8.5\p@{}%
			\right.%
			\n@space
			$%
		}%
	}%
}

\def\twelvebig@#1{%
	{%
		\hbox{%
			$%
			\left
			#1%
			\vbox to10.2\p@{}
			\right.%
			\n@space
			$%
		}%
	}%
}

%
%
%
%
%
\newif\ifl@beloutopen
\newwrite\l@belout
\newread\l@belin

\global\let\currentfile=\jobname

\def\getfile#1{%
	\immediate\closeout\l@belout
	\global\l@beloutopenfalse
	\gdef\currentfile{#1}%
	\input #1%
	\par
	\newpage
}

\def\getxrefs#1{%
	\bgroup
		\def\gobble##1{}
		\edef\list@{#1,}%
		\def\gr@boff##1,##2\end{
			\openin\l@belin=##1.xref
			\ifeof\l@belin
			\else
				\closein\l@belin
				\input ##1.xref
			\fi
			\def\list@{##2}%
			\ifx\list@\empty
				\let\next=\gobble
			\else
				\let\next=\gr@boff
			\fi
			\expandafter\next\list@\end
		}%
		\expandafter\gr@boff\list@\end
	\egroup
}

\def\testdefined#1#2#3{%
	\expandafter\ifx
	\csname #1\endcsname
	\relax
	#3%
	\else #2\fi
}

\def\document{%
	\minaw@11.11128\ex@ 
	\def\alloclist@{\empty}%
	\def\fontlist@{\empty}%
	\openin\l@belin=\jobname.xref	
	\ifeof\l@belin\else
		\closein\l@belin
		\input \jobname.xref
	\fi
}

\def\getst@te#1#2{%
	\edef\st@te{\csname #1s!#2\endcsname}%
	\expandafter\ifx\st@te\relax
		\def\st@te{0}%
	\fi
}

\def\setst@te#1#2#3{%
	\expandafter
	\gdef\csname #1s!#2\endcsname{#3}%
}

\outer\def\setupautolabel#1#2{%
	\def\newcount@{\global\alloc@0\count\countdef\insc@unt}	
	\def\newtoks@{\global\alloc@5\toks\toksdef\@cclvi}
	\expandafter\newcount@\csname #1Number\endcsname
	\expandafter\global\csname #1Number\endcsname=1%
	\expandafter\newtoks@\csname #1l@bel\endcsname
	\expandafter\global\csname #1l@bel\endcsname={#2}%
}

\def\reflabel#1#2{%
	\testdefined{#1l@bel}
	{
		\getst@te{#1}{#2}%
		\ifcase\st@te
			???
			\message{Unresolved forward reference to
				label #2. Use another pass.}%
		\or	
			\setst@te{#1}{#2}2
			\csname #1l!#2\endcsname 
		\or	
			\csname #1l!#2\endcsname 
		\or	
			\csname #1l!#2\endcsname 
		\fi
	}{
		{\escapechar=-1 
		\errmessage{You haven't done a
			\string\\setupautolabel\space for type #1!}%
		}%
	}%
}

{\catcode`\{=12 \catcode`\}=12
	\catcode`\[=1 \catcode`\]=2
	\xdef\Lbrace[{]
	\xdef\Rbrace[}]%
]%

\def\setlabel#1#2{%
	\testdefined{#1l@bel}
	{
		\edef\templ@bel@{\expandafter\the
			\csname #1l@bel\endcsname}%
		\def\@rgtwo{#2}%
		\ifx\@rgtwo\empty
		\else
			\ifl@beloutopen\else
				\immediate\openout\l@belout=\currentfile.xref
				\global\l@beloutopentrue
			\fi
			\getst@te{#1}{#2}%
			\ifcase\st@te
			\or	
			\or	
				\edef\oldnumber@{\csname #1l!#2\endcsname}%
				\edef\newnumber@{\templ@bel@}%
				\ifx\newnumber@\oldnumber@
				\else
					\message{A forward reference to label 
						#2 has been resolved
						incorrectly.  Use another
						pass.}%
				\fi
			\or	
				\errmessage{Same label #2 used in two
					\string\setlabel s!}%
			\fi
			\expandafter\xdef\csname #1l!#2\endcsname
				{\templ@bel@}
			\setst@te{#1}{#2}3%
			\immediate\write\l@belout 
				{\string\expandafter\string\gdef
				\string\csname\space #1l!#2%
				\string\endcsname
				\Lbrace\templ@bel@\Rbrace
				}%
			\immediate\write\l@belout 
				{\string\expandafter\string\gdef
				\string\csname\space #1s!#2%
				\string\endcsname
				\Lbrace 1\Rbrace
				}%
		\fi
		\templ@bel@	
		\expandafter\ifx\envir@end\endref 
			\gdef\marginalhook@{\marginal{#2}}%
		\else
			\marginal{#2}
		\fi
		\expandafter\global\expandafter\advance	
			\csname #1Number\endcsname
			by 1 %
	}{
		{\escapechar=-1
		\errmessage{You haven't done a \string\\setupautolabel\space
			for type #1!}%
		}%
	}%
}


\newcount\SectionNumber
\setupautolabel{t}{\number\SectionNumber.\number\tNumber}
\setupautolabel{r}{\number\rNumber}
\setupautolabel{T}{\number\TNumber}

\define\rref{\reflabel{r}}
\define\tref{\reflabel{t}}

\define\tnum{\setlabel{t}}
\define\rnum{\setlabel{r}}

%
\def\strutdepth{\dp\strutbox}%
\def\strutheight{\ht\strutbox}%

\newif\iftagmode
\tagmodefalse

\let\old@tagform@=\tagform@
\def\tagform@{\tagmodetrue\old@tagform@}

\def\marginal#1{%
	\ifvmode
	\else
		\strut
	\fi
	\ifdraftmode
		\ifmmode
			\ifinner
				\let\Vorvadjust=\Vadjust
			\else
				\let\Vorvadjust=\vadjust
			\fi
		\else
			\let\Vorvadjust=\Vadjust
		\fi
		\iftagmode	
			\llap{%
				\smalltype
				\vtop to 0pt{%
					\pretolerance=2000
					\tolerance=5000
					\raggedright
					\hsize=.72in
					\parindent=0pt
					\strut
					#1%
					\vss
				}%
				\kern.08in
				\iftagsleft@
				\else
					\kern\hsize
				\fi
			}%
		\else
			\Vorvadjust{%
				\kern-\strutdepth 
				{%
					\smalltype
					\kern-\strutheight 
					\llap{%
						\vtop to 0pt{%
							\kern0pt
							\pretolerance=2000
							\tolerance=5000
							\raggedright
							\hsize=.5in
							\parindent=0pt
							\strut
							#1%
							\vss
						}%
						\kern.08in
					}%
					\kern\strutheight
				}%
				\kern\strutdepth
			}
		\fi
	\fi
}


\newbox\Vadjustbox

\def\Vadjust#1{
	\global\setbox\Vadjustbox=\vbox{#1}%
	\ifmmode
		\ifinner
			\innerVadjust
		\fi		
	\else
		\innerVadjust
	\fi
}

\def\innerVadjust{%
	\def\nexti{\aftergroup\innerVadjust}%
	\def\nextii{%
		\ifvmode
			\hrule height 0pt 
			\box\Vadjustbox
		\else
			\vadjust{\box\Vadjustbox}%
		\fi
	}%
	\ifinner
		\let\next=\nexti
	\else
		\let\next=\nextii
	\fi
	\next
}%

\global\let\marginalhook@\empty

\def\endref{%
\setbox\tw@\box\thr@@
\makerefbox?\thr@@{\endgraf\egroup}%
  \endref@
  \endgraf
  \endgroup
  \keyhook@
  \marginalhook@
  \global\let\keyhook@\empty 
  \global\let\marginalhook@\empty 
}

\catcode`\@=\oldatsign


\topmatter
\title  Exact pairs of homogeneous zero divisors\endtitle 
\leftheadtext{Kustin, Striuli, and Vraciu}
\rightheadtext{Exact pairs of homogeneous zero divisors}
\author   
Andrew R. Kustin\footnote{Supported in part by the National Security Agency  and the Simons Foundation.\phantom
{xxxxxxxxxxxxxx}},
 Janet Striuli\footnote{Supported in part by the National Science Foundation. \phantom{xxxxxxxxxxxxxxxxxxxxxxxxxxxxxxxxxx}}, and Adela Vraciu\footnote{Supported in part by the National Security Agency and the National Science Foundation.\phantom
{xxxxxxx}}
\endauthor

\address
Mathematics Department,
University of South Carolina,
Columbia, SC 29208\endaddress
\email kustin\@math.sc.edu \endemail
\address Mathematics Department,
Fairfield University, 
Fairfield, Connecticut 06824 \endaddress
\email  jstriuli\@fairfield.edu\endemail
 \address
Mathematics Department,
University of South Carolina,
Columbia, SC 29208\endaddress
\email vraciu\@math.sc.edu \endemail
\keywords Compressed level algebra, Determinantal ring, Exact pair of zero divisors, Generic points in projective space, Linear resolution, Matrix factorization, Pfaffians, Segre embedding, Tate resolution, Totally acyclic complex, Totally reflexive module\endkeywords
\subjclass\nofrills{2010 {\it Mathematics Subject Classification.}}  13D02, 13A02
\endsubjclass
\abstract Let $S$ be a standard graded Artinian algebra over a field $k$. We identify constraints on the Hilbert function of $S$ which are imposed by the hypothesis that $S$ contains an exact pair of homogeneous  zero divisors.  As a consequence, we prove that if $S$ is a compressed level algebra, then $S$ does not contain any homogeneous  zero divisors.
\endabstract

\endtopmatter

\document
In \cite{\rref{HS}},  Henriques and  \c{S}ega   defined the pair of elements $(a,b)$ in a commutative ring $S$ to be an {\it exact pair of zero divisors} if $(0:_Sa)=(b)$ and $(0:_Sb)=(a)$. We take  $S$ be a standard graded Artinian algebra over a field and we identify constraints on the Hilbert function of $S$ which are imposed by the hypothesis that $S$ contains an exact pair  $(\theta_1,\theta_2)$  of homogeneous  zero divisors. In Theorem \tref{Conj} we prove that the main numerical constraint depends on the sum $\deg \theta_1+\deg \theta_2$, but not on the individual numbers $\deg \theta_1$ or $\deg \theta_2$. In other words, the numerical constraint imposed on $S$ by having an exact pair of homogeneous  zero divisors of degrees $d_1$ and $d_2$ is the same as the constraint imposed by having an exact pair of homogeneous  zero divisors of degrees $1$ and $d_1+d_2-1$. This result is especially curious because it is possible for $S$ to have an exact pair of homogeneous  zero divisors of degrees $2$ and $2$ without having any homogeneous exact zero divisors of degree $1$; see Example \tref{may4}. 
Our main  result is Theorem \tref{Conj}.
\proclaim{Theorem \tref{Conj}} 
Let $S$ be a standard graded Artinian $k$-algebra. Suppose that $(\theta_1,\theta_2)$  is an  exact pair of homogeneous zero divisors in $S$. If $D=\deg \theta_1+\deg \theta_2$, then the Hilbert series of $S$ is divisible by $\frac{t^D-1}{t-1}$. 
 \endproclaim

In the statement of Theorem \tref{Conj}, the algebra $S$ is Artinian, so the Hilbert series, $\operatorname{HS}_S(t)$, of $S$ is a polynomial in $\Bbb Z[t]$, the expression $\frac{t^D-1}{t-1}$ is equal to the polynomial $1+t+t^2+\dots+t^{D-1}$ of $\Bbb Z[t]$, and 

$$\matrix \format\l\\\text{``the Hilbert series of $S$ is divisible by $\frac{t^D-1}{t-1}$'' means that the polynomial}\\ \text{$1+t+t^2+\dots+t^{D-1}$ divides the polynomial $\operatorname{HS}_S(t)$ in the polynomial}\\\text{ring $\Bbb Z[t]$.}\endmatrix\tag\tnum{div}$$

We apply Theorem \tref{Conj} in Section 3 to obtain a list of conditions on the standard graded $k$-algebra $S$, each of which lead to the conclusion that $S$ does not have an exact pair of homogeneous zero divisors. These results are striking due to the connection between the existence of totally reflexive $S$-modules and the existence of exact zero divisors in $S$.   
A finitely generated $S$-module $M$ is called {\it totally reflexive} if there exists a doubly infinite sequence of
finitely generated free $S$-modules 
$$F:\quad \cdots \to F_1\to F_0 \to F_{-1}\to \cdots $$
such that M is isomorphic to the module $\operatorname{Coker}(F_1 \to F_0)$, and such that both $F$ and the dual sequence
$\operatorname{Hom}_S (F , S)$ are exact. The complex $F$ is called {\it totally acyclic}. For example, if $\theta_1$ and $\theta_2$ are a pair of exact zero divisors in $S$, then the complex 
$$F:\quad \cdots @>\theta_2>> S@>\theta_1>>  S@>\theta_2>>  S@>\theta_1>> \cdots $$is totally acyclic and the $S$-modules $S/(\theta_1)$ and $S/(\theta_2)$ are totally reflexive. 
Totally reflexive modules were first studied by Auslander and Bridger \cite{\rref{AB}}, who
proved that $S$ is Gorenstein if and only if every $S$-module has a totally reflexive syzygy. Over a
Gorenstein ring, the totally reflexive modules are precisely the maximal Cohen-Macaulay modules,
and these have been studied extensively. A main result
of Christensen, Piepmeyer, Striuli, and Takahashi \cite{\rref{CPST}, Thm. B} asserts that if $S$ is not
Gorenstein, then the existence of one non-free totally reflexive $S$-module implies the existence of
infinitely many non-isomorphic indecomposable totally reflexive $S$-modules. The proof in  \cite{\rref{CPST}} is not constructive; however many  methods   \cite{\rref{Y},\rref{H11},\rref{CJRSW},\rref{B12}} have been found for constructing non-isomorphic indecomposable totally reflexive $S$-modules. Most of these methods, especially those in \cite{\rref{H11}} and \cite{\rref{CJRSW}}, have  involved the use of a pair of exact zero divisors. Indeed, one result in \cite{\rref{CJRSW}} gives conditions on $S$ for which the existence of a non-free totally reflexive $S$-module implies the existence of an exact zero divisor in $S$. Furthermore, \cite{\rref{CJRSW}, Section 8} reformulates results of Conca \cite{\rref{C}}, Hochster and Laksov \cite{\rref{HL}}, and Yoshino \cite{\rref{Y}} to show that, under appropriate hypotheses, a generic standard graded algebra over an infinite field has an exact zero divisor.

We now describe the proof of Theorem \tref{Conj}. Let $d_i=\deg \theta_i$  and for integer $i$, define
$$\sigma_i=\sum_{j\equiv i \mod D} \operatorname{HF}_S(j),$$where  $\operatorname{HF}_S(\underline{\phantom{x}})$ is the Hilbert function of $S$. The indexed list $\{\sigma_i\}$ is periodic of period at most $D$. To prove the result we must   show that $\{\sigma_i\}$ has period $1$. 
The proof is carried out as follows.The existence of $(\theta_1,\theta_2)$ allows us to construct a totally acyclic complex $\Bbb F$. Each homogeneous strand of $\Bbb F$ is a finite exact sequence of vector spaces. The observation that vector space dimension is additive on finite exact sequences leads to the conclusion that  $\{\sigma_i\}$ has period at most $d_1$; and hence, $\{\sigma_i\}$ has period at most the $\operatorname{gcd}\{d_1,d_2\}$.  Henceforth, we assume that  $d_1$ and $d_2$ have a non-unit factor in common. 

One naive approach to proving Theorem \tref{Conj} would involve looking for a homogeneous factorization $\theta_1=\theta\cdot\check{\theta}$, with $\deg \theta=1$, and then considering the resulting exact pair of zero divisors $(\theta,\check{\theta}\cdot \theta_2)$. This approach is doomed to fail; however, if one looks for a matrix factorization of $\theta_1$, instead of a factorization of $\theta_1$ in $S$, then this naive approach does indeed work. In Lemma \tref{Tate}, we use the idea of the Tate resolution of the residue field of a hypersurface ring to create a matrix factorization $(M,\check{M})$ of $\theta_1\cdot\operatorname{id}$. At this point, we have a totally acyclic complex, periodic of period two, whose maps are 
$M\check{M}$ and  $\theta_2\cdot\operatorname{id}$. We prove in Lemma \tref{Lem2} that the maps $\check{M}$ and  $\theta_2\cdot M$ is also gives rise to a totally acyclic complex, periodic of period two. In Lemma \tref{May8} we combine Lemmas \tref{Tate} and \tref{Lem2} and give a recipe for 
using an exact pair of  zero divisors to 
build a numerically interesting totally acyclic complex $\Bbb G$.  By looking at graded strands of $\Bbb  G$, we obtain the equations
$$\sum _{\ell=0}^{s_1} (-1)^{\ell}\binom{s_1}{\ell}\sigma_{N-\ell}=0$$ for all integers $N$, where $s_1=\operatorname{HS}_S(1)$. The coefficient matrix for this system of linear equations is a circulant matrix; in Lemma \tref{May9} we show that $\{\sigma_i\}$ has period one is the only solution.  

All of our notation and conventions are explained in Section 1. The details in the proof of Theorem \tref{Conj} are given immediately after the statement at the end of Section 2. Section 3 consists of examples and applications of Theorem \tref{Conj}. In Proposition \tref{MIF*} we show that, in general, compressed level algebras do not have any homogeneous exact zero divisors. Assorted examples of compressed level algebras are given. These examples include Gorenstein rings with linear resolutions, certain rings defined by Pfaffians, certain determinantal rings, and certain rings arising from  sets of generic points in projective space.  We also exhibit families of 
  standard graded Artinian $k$-algebras, which are not compressed level algebras, but which nonetheless  do not contain any homogeneous exact zero divisors. These families include rings which arise from Segre embeddings and more  determinantal rings.

\SectionNumber=1\tNumber=1
\heading Section \number\SectionNumber. \quad Terminology, notation, and preliminary results.
\endheading

We use $\operatorname{gcd}$ as an abbreviation for {\it greatest common divisor},  $\Bbb Z$ to represent the set of integers $\{\dots,-2,-1,0,1,2,\dots\}$,  $\Bbb N$ the set of positive integers $\{1,2,3,\dots\}$, and $\Bbb Q$ is the field of rational numbers. 

If $S$ is a ring, $N\subseteq M$ are $S$-modules and $X$ is a subset of $M$, then $$N:_SX=\{s\in S\mid sx\in N\text{ for all $x\in X$}\}.$$

\definition{Conventions \tnum{D1}}Let $k$ be a field. If $V$ is a vector space over $k$, then $\dim_kV$ is the vector space dimension of $V$. 

\smallskip\item {(1)}A {\it standard graded  $k$-algebra} is a graded ring $S=\bigoplus_{i\in\Bbb Z}[S]_i$, with $[S]_i=0$ for $i<0$, $[S]_0=k$, $\dim_k[S]_1<\infty$, and $S$ is generated as a $k$-algebra by $[S]_1$. 

\smallskip\item{(2)} For each   graded $S$-module $M$ we use $[M]_i$ to denote the {\it homogeneous component of $M$ of degree $i$}. 

\smallskip\item{(3)} The  {\it Hilbert function} of $M$ is the function $\operatorname{HF}_S(M,\underline{\phantom{x}})$, from the set $\Bbb Z$ to the set $\{0\}\cup \Bbb N\cup \{\infty\}$, with
$\operatorname{HF}_S(M,i)=\dim_k[M]_i$. We abbreviate $\operatorname{HF}_S(S,i)$ as  $\operatorname{HF}_S(i)$. 

\smallskip\item{(4)} The {\it Hilbert series} of a graded, finitely generated, $S$-module $M$ is the formal generating function $\operatorname{HS}_S(M,t)=\sum_{i\in \Bbb Z}\operatorname{HF}_S(M,i)t^i$. We abbreviate $\operatorname{HS}_S(S,t)$ as $\operatorname{HS}_S(t)$; this formal power series is called the {\it Hilbert series of $S$}. {\bf If $S$ is Artinian, then the ``Hilbert series'' of $S$ is actually a polynomial.}

\smallskip\item{(5)} If $M$ is a graded $S$-module and $a$ is an integer, then $M(a)$ is the graded $S$-module with $[M(a)]_i=[M]_{a+i}$ for all integers $i$.

\smallskip\item{(6)} If $R$ is a  standard graded $k$-algebra of dimension $d$, and $\ell_1,\dots,\ell_d$ is a regular sequence of linear forms in $R$, then    $S=R/(\ell_1,\dots,\ell_d)$ is called an {\it Artinian reduction} of $R$.\enddefinition

\definition{Definition \tnum{D1.5}} A complex of modules 
$$(\Bbb F,f)\: \cdots @> f_2>> \Bbb F_{1}@> f_1>> \Bbb  F_{0} @> f_0>> \Bbb F_{-1} @> f_{-1}>> \cdots$$ over the ring $S$ is {\it acyclic} if the homology $\operatorname{H}(\Bbb F)$ is equal to zero. The acyclic complex $\Bbb F$ is {\it totally acyclic} if $\Bbb F^*=\operatorname{Hom}_S(\Bbb F,S)$ is also acyclic. \enddefinition

\definition{Conventions \tnum{D2}} Let $S$ be a standard graded algebra over a field $k$ and $\varphi\:F\to G$ be a homomorphism of finitely generated free  graded $S$-modules. The homomorphism $\varphi$ is called {\it homogeneous of degree $d$} if, whenever $x$ is a homogeneous element of $F$, then $\varphi(x)$ is a homogeneous element of $G$ of degree $d+\deg x$. If $\varphi$ is called a {\it homogeneous homomorphism} and no degree is specified, then $\varphi$ is a homogeneous homomorphism of degree $0$. The homogenous homomorphism $\varphi\:F\to G$ is a {\it minimal homogeneous homomorphism} if $\varphi(F)\subseteq S_+G$,   where $S_+$ is the ideal $\sum_{1\le i}[S]_i$ of $S$.\enddefinition 

Observation \tref{may11} is our tool for extracting numerical information about the Hilbert function of a a standard graded Artinian algebra $S$ from the twists in a minimal homogeneous  totally acyclic complex of finitely generated free $S$-modules. The interesting part of the assertion is that the sum (\tref{sum}) is finite.
\proclaim{Observation \tnum{may11}} Let $S$ be a standard graded Artinian algebra over the field $k$ and $(\Bbb F,f)$ be a minimal homogeneous  totally acyclic complex of finitely generated free $S$-modules. Write $\Bbb F$  in the form 
$$\dots@>f_{j+1}>> \Bbb F_j@> f_j >> \Bbb F_{j-1}@> f_{j-1} >>\dots\ .$$ Then, for each integer $N$, the expression
$$\sum _{p\in \Bbb Z} \left(\operatorname{HF}_S(\Bbb F_{2p},N)-\operatorname{HF}_S(\Bbb F_{2p-1},N)\vphantom{E^E}\right)\tag\tnum{sum}$$
is a finite sum and is equal to zero.\endproclaim

\demo{Proof} The free $S$-module $\Bbb F_j$ occupies the homological position $j$ in the complex $\Bbb F$, where $j$  varies over all of $\Bbb Z$. For each integer $N$, the homogeneous component $[\Bbb F]_N$ of $\Bbb F$ of degree $N$ is $\bigoplus_{j\in \Bbb Z}[\Bbb F_j]_N$; this component is an exact sequence of finite dimensional vector spaces
$$ \dots \to  [\Bbb F_{j+1}]_N\to [\Bbb F_j]_N \to [\Bbb F_{j-1}]_N\to \dots\ . $$We show below that $$\text{once $N$ is fixed, then $[\Bbb F_j]_N$ is non-zero for only a finite number of $j$.}\tag\tnum{below}$$ Vector space dimension is additive on finite exact sequences; thus, for each integer $N$, 
$$\split 0= \sum_{j\in \Bbb Z} (-1)^j\dim_k[\Bbb F_j]_N&{}=
\sum_{p\in \Bbb Z}(\dim_k[\Bbb F_{2p}]_N-\dim_k[\Bbb F_{2p-1}]_N)\\&{}=\sum _{p\in \Bbb Z} \left(\operatorname{HF}_S(\Bbb F_{2p},N)-\operatorname{HF}_S(\Bbb F_{2p-1},N)\vphantom{E^E}\right).\endsplit $$

We complete the argument by establishing (\tref{below}). The graded ring $S$ is Artinian; so $[S]_i$ is zero for all large $i$. Define $e$ to be the largest integer with $[S]_e\neq 0$.  Fix an integer $j$. The free $S$-module $\Bbb F_j$ is graded and finitely generated. Let $r_j$ be the rank of $\Bbb F_j$ and $m_{j,1}\le \dots\le m_{j,r_j}$ be the degrees of the elements in a homogeneous minimal generating set for $\Bbb F_j$. We are particularly interested in $\ell_j=m_{j,1}$ and $h_j=m_{j,r_j}$. Notice that
$$[\Bbb F_j]_N=0\quad\text{for all $N$ with either $N<\ell_j$ or  $h_j+e<N$}.\tag\tnum{lohi}$$
The key to establishing   (\tref{below})  is contained in the inequalities $$\ell_j+1\le \ell_{j+1}\quad\text{and}\quad h_{j-1}\le h_j-1.\tag\tnum{next}$$ Roughly speaking, the inequality on the left says that ``no column of $f_{j+1}$ can consist entirely of zeros'' and the inequality on right says that ``no row of $f_{j-1}$ can consist entirely of zeros''.  
To show the equality on the left side of (\tref{next}), we consider 
 a homogeneous minimal generator   $\xi$ of $\Bbb F_{j+1}$. The complex $\Bbb F$ is minimal,  homogeneous, and acyclic; thus, $f_{j+1}(\xi)$ is not zero  and   is not a minimal generator of $\Bbb F_j$. 
Apply (\tref{lohi}) to conclude, in the first place,  that   $\ell_j\le \deg f_{j+1}(\xi)$. Every homogeneous element in $\Bbb F_j$ of degree $\ell_j$ is a minimal generator of $\Bbb F_j$; so, we also conclude that  $\ell_j\neq \deg f_{j+1}(\xi)$; therefore, we have shown that $\ell_j+1\le \deg f_{j+1}(\xi)$. The elements $\xi$ and $f_{j+1}(\xi)$ of $\Bbb F$ have the same degree and the equality on the left side of (\tref{next}) is established. To prove the inequality on the right side of (\tref{next}), we consider the dual complex $$\Bbb F^*:\quad \dots @> f_{j-1}^*>> \Bbb F^*_{j-1} @> f_{j}^*>> \Bbb F^*_{j} @> f_{j+1}^*>>\dots,$$ which is also homogeneous, minimal, and acyclic. The elements of a minimal homogeneous generating set  for $\Bbb F_j^*$ have degrees 
$$-h_j=-m_{j,r_j}\le \dots\le -m_{j,1}\le -\ell_j.$$ Apply the inequality on the left side of (\tref{next}) to see that $-h_j+1\le -h_{j-1}$, which is the inequality on the right side of (\tref{next}).  

Suppose that $[\Bbb F_j]_N\neq 0$. It follows from (\tref{lohi}) that $\ell_j\le N\le h_j+e$. In particular, the integers $N-\ell_j$ and $h_j+e-N$ are non-negative. If $b$ is an integer with $N-\ell_j<b$, then (\tref{next}) gives $N<\ell_j+b\le \ell_{j+b}$; hence, (\tref{lohi}) yields $[\Bbb F_{j+b}]_N=0$. Similarly, if $h_j+e-N<b$, then (\tref{next}) gives $h_{j-b}+e\le h_j-b+e<N$ and (\tref{lohi}) yields $[\Bbb F_{j-b}]_N=0$. \qed \enddemo

\SectionNumber=2\tNumber=1
\heading Section \number\SectionNumber. \quad The proof of the main Theorem.
\endheading 

We first produce a matrix factorization of an essentially arbitrary element of a ring. Our technique is inspired by the Tate resolution of the residue field of a hypersurface ring. 

\proclaim{Lemma \tnum{Tate}} 
 Let be $S$ a standard graded algebra over a field $k$ and 
  $\theta$ be  a homogeneous element of $S$ of degree $d\ge 2$. Then there exist finitely generated, free, graded $S$-modules $F$ and $G$ and minimal homogeneous $S$-module homomorphisms $M\:F\to G$ and $\check{M}\:G\to F(d)$ such that the compositions $\check{M}\circ M\:F\to F(d)$ and $M\circ \check{M}\: G\to G(d)$ both are multiplication by $\theta$. \endproclaim

\demo{Proof} Let $x_1,\dots,x_{s_1}$ be a basis for the vector space $[S]_1$. Identify elements $y_1,\dots,y_{s_1}$ in  $[S]_{d-1}$ with $\theta=\sum_{i=1}^{s_1}x_iy_i$. Let $V$ be the graded free $S$-module $\bigoplus_{i=1}^{s_1} S\varepsilon_i$, where each $\varepsilon_i$ has degree $1$, $\Theta$ be a divided power variable of degree $2$, and $\Bbb T$ be the  Graded Divided-power   Algebra $\Bbb T=(\bigwedge_S^{\bullet}V){<}\Theta{>}$. Define a map $t\:\Bbb T\to \Bbb T$ as follows:  $t(\varepsilon_i)=x_i\in S$, the restriction of $t$ to $\bigwedge_S^{\bullet}V$ is the usual Koszul complex map associated to $t\:V\to S$, and if $\kappa$ is in $\bigwedge_S^iV$ and $1\le\ell$, then
$$t(\kappa\otimes \Theta^{(\ell)})=t(\kappa)\otimes \Theta^{(\ell)} + (-1)^i \kappa\wedge ({\tsize\sum_j}y_j\varepsilon_j)\otimes\Theta^{(\ell-1)}.$$ Notice that $t(t(\Theta))=\theta$. Fix an integer  $p$  with $s_1+2\le 2p$, then
$$\Bbb T_{2p}=\bigoplus_{i=0}^{\lfloor\frac {s_1}2\rfloor} {\tsize\bigwedge_S^{2i}}V\otimes_SS\Theta^{(p-i)}, \quad 
\Bbb T_{2p-1}=\bigoplus_{i=0}^{\lfloor\frac {s_1-1}2\rfloor} {\tsize\bigwedge_S^{2i+1}}V\otimes_SS\Theta^{(p-1-i)},
\quad\text{and}$$
$$(t\circ t)(\kappa\otimes \Theta^{(\ell)})=\theta\cdot (\kappa\otimes \Theta^{(\ell-1)})\quad\text{for all $\kappa\otimes \Theta^{(\ell)}$ in $\Bbb T_{2p}$ or $\Bbb T_{2p-1}$.}$$ We see, from the definition of $t$, that $t\:V\to S$ may be written as  $S(-1)^{s_1}\to S$ and $t\: S\Theta\to V$ may be written as $S(-d)^1\to S(-1)^{s_1}$.   Thus, 
$$\allowdisplaybreaks\align\Bbb T_{2p}&{}=\bigoplus_{i=0}^{\lfloor\frac {s_1}2\rfloor}S(-2i-(p-i)d)^{\binom {s_1}{2i}}\tag{\tnum{?}}\\&{}=
\bigoplus_{i=0}^{\lfloor\frac {s_1}2\rfloor}S(i(d-2)-pd)^{\binom {s_1}{2i}}
\quad\text{and}\\
\Bbb T_{2p-1}&{}=\bigoplus_{i=0}^{\lfloor\frac {s_1-1}2\rfloor}S(-2i-1-(p-1-i)d)^{\binom {s_1}{2i+1}}\\
&{}=\bigoplus_{i=0}^{\lfloor\frac {s_1-1}2\rfloor}S(i(d-2)+d-1-pd)^{\binom {s_1}{2i+1}}.\endalign$$
Notice that, in the language of (\tref{?}), $\Bbb T_{2p-2}=\Bbb T_{2p}(d)$, $\Bbb T_{2p-3}=\Bbb T_{2p-1}(d)$, and  $t\:\Bbb T_{2p-2}\to \Bbb T_{2p-3}$ is a shift of $t\:\Bbb T_{2p}\to \Bbb T_{2p-1}$. Define $$F= \Bbb T_{2p}(pd), \quad G=\Bbb T_{2p-1}(pd),$$ $M$ to be $t\:F\to G$, and $\check{M}$ to be $t\:G\to F(d)$. We have finished the proof. For future reference we record 
$$F=
\bigoplus_{i=0}^{\lfloor\frac {s_1}2\rfloor}S(i(d-2))^{\binom {s_1}{2i}}
\quad\text{and}\quad
G=\bigoplus_{i=0}^{\lfloor\frac {s_1-1}2\rfloor}S(i(d-2)+d-1)^{\binom {s_1}{2i+1}}. \qed\tag\tnum{conl}$$ 
\enddemo

Suppose that $M$, $\check{M}$, and $M_2$ are three maps, which satisfy sufficient commuting relations, and for which
there is a totally acyclic complex, periodic of period two, whose maps are 
$M\check{M}$ and  $M_2$. We prove in Lemma \tref{Lem2} that the maps $\check{M}$ and  $M_2  M$ is also gives rise to a totally acyclic complex, periodic of period two.

\proclaim{Lemma \tnum{Lem2}} Let be $S$ a standard graded algebra over a field $k$, $d_1$ and $d_2$ be positive integers, $F$ and $G$ be finitely generated, graded, free $S$-modules, and $M\:F\to G$, $\check{M}\: G\to F(d_1)$  be homogeneous $S$-module homomorphisms. Suppose further that 
  $M_2$ is a  matrix with entries from $S$ such that  $M_2\:F\to F(d_2)$ and $M_2\:G\to G(d_2)$ both are homogeneous $S$-module homomorphisms. Assume that the diagram $$\CD F@> M>> G@> \check{M}>> F(d_1)\\@V M_2VV @V M_2 VV @V M_2VV \\F(d_2)@> M>> G(d_2)@> \check{M}>> F(d_1+d_2)\endCD$$ commutes and that the homomorphisms
$$\eightpoint \split&\Bbb F\:\quad \cdots @> \check{M}M >>F(-d_2) @> M_2 >> F@> \check{M}M >> F(d_1)@> M_2>>F(d_1+d_2) @> \check{M}M >> F(2d_1+d_2) @> M_2>> \cdots \ \text{and}\\&\check{\Bbb F}\:\quad \cdots @> M\check{M} >>G(-d_2) @> M_2>>  G@> M\check{M} >> G(d_1)@> M_2>>G(d_1+d_2) @> M\check{M} >> G(2d_1+d_2)@> M_2>> \cdots\phantom{\ \text{and}}\endsplit $$
form acyclic complexes. Then the homomorphisms
$$\eightpoint \split&\Bbb G\:\quad \cdots @> MM_2 >>G(- d_1) @> \check{M} >> F@> MM_2 >> G(d_2)@> \check{M}>>F(d_1+d_2) @> MM_2 >> G(d_1+2d_2) @> \check{M}>> \cdots \ \text{and}\\&\check{\Bbb G}\:\quad\cdots @> \check{M}M_2 >>F @> M>>  G@> \check{M}M_2 >> F(d_1+d_2)@> M>>G(d_1+d_2) @> \check{M}M_2 >> F(2d_1+2d_2)@> M>> \cdots\phantom{\ \text{and}}\endsplit $$
form acyclic complexes. Furthermore, if $\Bbb F$ and $\check{\Bbb F}$ are totally acyclic complexes, then $\Bbb G$ and $\check{\Bbb G}$ also are totally acyclic complexes.
\endproclaim

\demo{Proof} 

\flushpar{\bf Claim 1.} $\operatorname{Ker} M=\operatorname{Im}(\check{M}M_2)$

The inclusion $\supseteq$ is obvious. We show $\subseteq$. Take $x\in F(a)$, for some integer $a$, with $M(x)=0$. It follows that $x$ is in $\operatorname{Ker} (\check{M}M)=\operatorname{Im} M_2$. Thus, there is an element $x_1$ of $F(a-d_2)$, with $x=M_2x_1$. The hypothesis $x\in \operatorname{Ker} M$ now gives
$$0=Mx=MM_2x_1.$$The maps $M$ and $M_2$ commute; so $Mx_1\in \operatorname{Ker} M_2=\operatorname{Im}(M\check{M})$; and $$Mx_1=M\check{M}x_2$$ for some $x_2$ in $G(a-d_1-d_2)$. We see that $x_1=\check{M}x_2+x_3$, where $x_3$ is the   element $x_1-\check{M}x_2$ of $\operatorname{Ker} M$. We have
$$x=M_2 x_1=M_2(\check{M}x_2+x_3)\in \operatorname{Im} (M_2\check{M})+M_2(\operatorname{Ker} M).$$ In other words, $\operatorname{Ker} M\subseteq \operatorname{Im} (M_2\check{M})+M_2 \operatorname{Ker} (M)$. Iterate the above argument to see that $\operatorname{Ker} M\subseteq \operatorname{Im} (M_2\check{M})+M_2^r \operatorname{Ker} (M)$ for all positive integers $r$. On the other hand, the maps all are homogeneous; so, for each fixed $i$, $$[\operatorname{Ker} M]_i\subseteq [\operatorname{Im} (M_2\check{M})]_i+M_2^r ([\operatorname{Ker} (M)]_{i-rd_2}).$$ The module $F$ is finitely generated; thus, when $r$ is sufficiently large, $[F]_{i-rd_2}=0$. It follows that $[\operatorname{Ker} M]_i\subseteq [\operatorname{Im} (M_2\check{M})]_i$, for all $i$; hence, $\operatorname{Ker} M\subseteq \operatorname{Im} (M_2\check{M})$ and this completes the proof of Claim 1.

\medskip\flushpar{\bf Claim 2.} $\operatorname{Ker} \check{M}=\operatorname{Im}(MM_2)$

One repeats the proof of Claim 1 after reversing the roles of $M$ and $\check{M}$.

\medskip\flushpar{\bf Claim 3.} $\operatorname{Ker} \check{M}M_2=\operatorname{Im}(M)$

The inclusion $\supseteq$ is obvious. We show $\subseteq$. Take $x\in G(a)$, for some integer $a$, with $x\in\operatorname{Ker} \check{M}M_2$. The hypothesis ensures that $M_2\check{M}x=\check{M}M_2x=0$. It follows that   $\check{M}x$ is in $\operatorname{Ker} M_2=\operatorname{Im}(\check{M}M)$ and $\check{M}x=\check{M}Mx_1$ for some $x_1\in F(a)$. Therefore, $x-M x_1\in \operatorname{Ker} \check{M}$ which we saw in Claim 2 is equal to $\operatorname{Im}(MM_2)$. We have $x-Mx_1=MM_2x_2$ for some $x_2\in F(a-d_2)$ and $x\in \operatorname{Im}(M)$. 

The  equality $\operatorname{Ker} MM_2=\operatorname{Im}(\check{M})$  follows from the symmetry of the situation the same way that Claim 2 followed from Claim 1. We conclude that $\Bbb G$ and $\check{\Bbb G}$ both are acyclic complexes. 

Suppose now that $\Bbb F$ and $\check{\Bbb F}$ are totally acyclic complexes. In this case, $\Bbb F^*$ and $\check{\Bbb F}^*$ are acyclic complexes, the maps $M^*$ and $M_2^*$ commute, and the maps $\check{M}^*$ and $M_2^*$ commute. One may apply what we have already shown to conclude that the complexes $\Bbb G^*$ and $\check{\Bbb G}^*$ are acyclic; and therefore,  $\Bbb G$ and $\check{\Bbb G}$ both are totally acyclic complexes. \qed
\enddemo

  In Lemma \tref{May8} we combine Lemmas \tref{Tate} and \tref{Lem2} and give a recipe for 
using an exact pair of  zero divisors to 
build a numerically interesting totally acyclic complex $\Bbb G$. 

\proclaim{Lemma \tnum{May8}} Let $S$ be a standard graded Artinian algebra over the field $k$. Suppose that $(\theta_1,\theta_2)$  is an exact  pair of  homogeneous zero divisors in $S$, with  $d_1=\deg \theta_1\ge 2$ and $d_2= \deg \theta_2$. Let $D=d_1+d_2$ and $s_1=\operatorname{HF}_S(1)$. Then there is a homogeneous totally acyclic complex $\Bbb G$ of free $S$-modules of the form
$$\dots @>\varphi>> \Bbb G_{2p+1}@> \psi >> \Bbb G_{2p}@>\varphi>> \Bbb G_{2p-1}@> \psi >> \dots, \quad\text{with}$$
$$\split \Bbb G_{2p}= {}&\bigoplus_{i=0}^{\lfloor\frac {s_1}2\rfloor}S(i(d_1-2)-pD)^{\binom {s_1}{2i}} \quad\text{and}\\
\Bbb G_{2p-1}= {}&\bigoplus_{i=0}^{\lfloor\frac {s_1-1}2\rfloor}S(i(d_1-2)+d_1-1-pD)^{\binom {s_1}{2i+1}},\endsplit
$$ for all integers $p$. The matrices $\varphi$ and $\psi$ are square with $2^{s_1-1}$ rows and columns. 
  \endproclaim

\demo{Proof}   Apply Lemma \tref{Tate} to the homogeneous element $\theta_1\in S$ of degree $d_1\ge 2$ to obtain 
finitely generated, free, graded $S$-modules $F$ and $G$ and minimal homogeneous $S$-module homomorphisms $M\:F\to G$ and $\check{M}\:G\to F(d_1)$ such that the compositions $\check{M}\circ M\:F\to F(d_1)$ and $M\circ \check{M}\: G\to G(d_1)$ both are multiplication by $\theta_1$. Let $M_2$ be the $S$-module homomorphism ``multiplication by $\theta_2$''. Notice that $M_2$ is a legitimate map from $F$ to $F(d_2)$ and $M_2$ is also a legitimate map from $G$ to $G(d_2)$. The hypotheses of Lemma \tref{Lem2} all are satisfied because  $(\theta_1,\theta_2)$  is an exact pair of zero divisors in $S$ and  the composition of $M$ and $\check{M}$, in either order, is multiplication by $\theta_1$.   We record the acyclic complex $\Bbb G$, from the conclusion of Lemma \tref{Lem2}, with 
$$\cdots@>\varphi>>\Bbb G_1=G(-d_1)@>\psi >>F=\Bbb G_0@> \varphi>> G(d_2)=\Bbb G_{-1}@>\psi >> \cdots \ ,$$  $\psi=\check{M}$, $\varphi=MM_2$, and, according to  (\tref{conl}),  
$$\split \Bbb G_{2p}&{}=F(-pD)=
\bigoplus_{i=0}^{\lfloor\frac {s_1}2\rfloor}S(i(d_1-2)-pD)^{\binom {s_1}{2i}}
\quad\text{and}\\
\Bbb G_{2p-1}&{}=G(-pD)=\bigoplus_{i=0}^{\lfloor\frac {s_1-1}2\rfloor}S(i(d_1-2)+d_1-1-pD)^{\binom {s_1}{2i+1}}. \qed\endsplit$$ 
\enddemo

Lemma \tref{May9} is a straightforward numerical calculation. The assertions about circulant matrices are easy to verify and are stated on   Wikipedia. It is gratifying to solve a system of linear equations without making any unpleasant calculations.

\proclaim{Lemma \tnum{May9}} Assume that 
\roster\item   $\sigma_i$ is an integer defined for all $i\in \Bbb Z$, \item $a$ is a fixed  positive integer with $\sigma_i=\sigma_{i+a}$ for all $i\in \Bbb Z$, \item and $b$  is a fixed positive integer with $\sum\limits_{i=0}^b(-1)^i\binom bi\sigma_{N-i}=0$ for all $N\in \Bbb Z$.\endroster Then $\sigma_i=\sigma_j$ for all integers $i$ and $j$.\endproclaim
\demo{Proof} Pick an integer $B$ such that  $b+1\le B$ and $B$ is  a multiple of $a$. Hypothesis (2) gives $$\sigma_{i+B}=\sigma_i \text{ for all $i\in \Bbb Z$}.\tag\tnum{B}$$ Define integers $c_0,\dots,c_{B-1}$ by 
$$c_i=\cases (-1)^i\binom bi&\text{if $0\le i\le b$}\\0&\text{if $b+1\le i\le B-1$},\endcases$$ 
and let $C$ be the $B\times B$ circulant matrix
$$C=\bmatrix c_0&c_1&c_2&\dots&c_{B-2}&c_{B-1}\\c_{B-1}&c_0&c_1&\dots&c_{B-3}&c_{B-2}\\
c_{B-2}&c_{B-1}&c_0&\dots&c_{B-4}&c_{B-3}\\&&\ddots&\ddots&\ddots\\
c_2&c_3&c_4&\dots&c_0&c_1 \\c_1&c_2&c_3&\dots&c_{B-1}&c_0\endbmatrix$$ and for each integer $N$, let $\Sigma_N$ be the column vector
$$\Sigma_N=\bmatrix \sigma_N\\\vdots\\\sigma_{N-B+1} \endbmatrix$$ in $\Bbb Z^B$. We may read hypothesis (3) to say
$$\text{(row 1 of $C$)} \Sigma_N=0,\text{ for all $N\in \Bbb Z$}.\tag\tnum{c}$$
Equation (\tref{c}) holds when $N$ is replaced by $N+1$ and formula (\tref{B}) yields that $\sigma_{N+1}=\sigma_{N-B+1}$; so
$\text{(row 2 of $C$)} \Sigma_N$ is also zero. One may iterate this procedure to see that $C\Sigma_N$=0. We complete the proof by showing that the kernel of $C$ is generated by the $B\times 1$  column vector with every entry equal to $1$. 

It is well known that $C$ is a diagonalizable matrix over the complex numbers. Let 
$\omega$ be a primitive $B^{\text{\rm th}}$ root of $1$, $p(x)$ be the polynomial $p(x)=\sum_{i=0}^{B-1}c_ix^i\in \Bbb Z[x]$, 
$V$ be the Vandermonde matrix 
matrix $$V=\bmatrix 1&1&1&\dots &1\\
1&\omega&\omega^2&\dots&\omega^{B-1}\\
\vdots\\
1&\omega^{B-1}&\omega^{2(B-1)}&\dots&\omega^{(B-1)(B-1)}\endbmatrix,$$ and $\Cal D$ be the diagonal matrix with diagonal entries
$p(1),p(\omega),\dots,p(\omega^{B-1})$. It is clear that $CV=V\Cal D$ and $V$ is invertible; so, the matrices $C$ and $\Cal D$ have the same rank. On the other hand, our choice of the coefficients $c_i$ gives $p(x)=(1-x)^b$. It follows that exactly one of the diagonal entries of $\Cal D$ is zero and the corresponding eigenvector is the first column of $V$, as we desired. \qed\enddemo

Theorem \tref{Conj} is the main result of this paper. 

\proclaim{Theorem \tnum{Conj}} 
Let $S$ be a standard graded Artinian algebra over a field $k$.
If $(\theta_1,\theta_2)$  is an  exact pair of homogeneous zero divisors in $S$ and $D=\deg \theta_1+\deg \theta_2$, then the Hilbert series of $S$ is divisible by $\frac{t^D-1}{t-1}$ in the sense of {\rm(\tref{div})}. \endproclaim

\remark{\bf Remark \tnum{alt}} The set of polynomials in $\Bbb Z[t]$ which are divisible by $$\frac{t^D-1}{t-1}=\sum_{i=0}^{D-1}t^i$$ forms an ideal $J$ and the polynomial $p(t)=\sum c_it^i$ in $\Bbb Z[t]$ is in $J$ if and only if 
$$\sum_{\ell\in \Bbb Z}c_{0+\ell D}=\sum_{\ell\in \Bbb Z}c_{1+\ell D}=\dots = \sum_{\ell\in \Bbb Z}c_{D-1+\ell D}. \tag\tnum{sigma}$$ View $c_i$ to be zero if $i$ is negative or greater than the degree of $p(t)$. The equations (\tref{sigma}) often are the easiest way to test if the Hilbert series of a given Artinian algebra satisfy the conditions of Theorem \tref{Conj}. When we apply (\tref{sigma}), we let $\sigma_i$ denote
$$\sigma_i=\sum_{\ell\in \Bbb Z}c_{i+\ell D}.$$ \endremark

\demo{Proof of Theorem {\rm\tref{Conj}}} Let $d_1=\deg \theta_1$ and $d_2=\deg \theta_2$. For each integer $i$, let $s_i=\operatorname{HF}_S(i)$. (The Hilbert function $\operatorname{HF}_S$ and the Hilbert series $\operatorname{HS}_S$ are both defined in Conventions \tref{D1}.) Define 
$$\sigma_i=\sum_{\ell\in \Bbb Z}s_{i+\ell D},\tag\tnum{sig}$$
for all integers $i$. The ring $S$ is Artinian; so, each $\sigma_i$ is a finite integer. The definition of $\sigma_i$ shows that $\sigma_i=\sigma_{i+D}$ for all integers $i$. In light of Remark \tref{alt}, it suffices to prove that 
$$\sigma_i=\sigma_j \quad\text{for all integers $i$ and $j$.}\tag\tnum{mini}$$
The hypothesis that $(\theta_1,\theta_2)$  is an exact  pair of  homogeneous zero divisors in $S$ ensures that there is a homogeneous totally acyclic complex $\Bbb F$ of free $S$-modules of the form
$$\dots @>\theta_1>> \Bbb F_{2p+1}@> \theta_2 >> \Bbb F_{2p}@>\theta_1>> \Bbb F_{2p-1}@> \theta_2 >> \dots,$$
with$$\Bbb F_{2p}=S(-pD)\quad\text{and}\quad \Bbb F_{2p-1}=S(-pD+d_1).$$ Apply Observation \tref{may11} to see that  
$$\split 0&{}=\sum _{p\in \Bbb Z} \left(\operatorname{HF}_S(\Bbb F_{2p},N)-\operatorname{HF}_S(\Bbb F_{2p-1},N)\vphantom{E^E}\right)
=\sum _{p\in \Bbb Z} \left(s_{N-pD}-s_{N-pD+d_1}\vphantom{E^E}\right)\\&{}=\sigma_N-\sigma_{N+d_1},\endsplit$$ for each integer $N$.
 It follows by symmetry that $\sigma_{N+d_2}=\sigma_{N+d_2-D}=\sigma_{N-d_1}=\sigma_N$; hence,
$$\sigma_{N+d_2}=\sigma_N=\sigma_{N+d_1}\quad\text{for all $N$.}\tag\tnum{i+d}$$
If the greatest common divisor of $d_1$ and $d_2$ is $1$, then there is nothing more to prove. Henceforth, we may assume that $d_1$ and $d_2$ have a non-unit factor in common. In particular, we may assume that $d_1\ge 2$. 
The fact that $S$ contains an exact  pair of homogeneous zero divisors of degrees $d_1$ and $d_2$, with $d_1\ge 2$, allows us to create the homogeneous totally acyclic complex $\Bbb G$ of Lemma \tref{May8}. Apply Observation \tref{may11}, together with (\tref{sig}) and (\tref{i+d}),   to see that 
$$\align 0&{}=\sum _{p\in \Bbb Z} \left(\operatorname{HF}_S(\Bbb G_{2p},N)-\operatorname{HF}_S(\Bbb G_{2p-1},N)\vphantom{E^E}\right)\\\allowdisplaybreak
&{}=\sum _{p\in \Bbb Z} \left(\sum_{i=0}^{\lfloor \frac{s_1}2\rfloor}\binom{s_1}{2i}s_{N+i(d_1-2)-pD}-\sum_{i=0}^{\lfloor \frac{s_1-1}2\rfloor}\binom{s_1}{2i+1} s_{N+i(d_1-2)+d_1-1-pD}\right)\\\allowdisplaybreak&{}
=\sum_{i=0}^{\lfloor \frac{s_1}2\rfloor}\binom{s_1}{2i}\left(\sum _{p\in \Bbb Z} s_{N+i(d_1-2)-pD}\right)-\hskip-1.81pt\sum_{i=0}^{\lfloor \frac{s_1-1}2\rfloor}\binom{s_1}{2i+1} \left(\sum _{p\in \Bbb Z} s_{N+i(d_1-2)+d_1-1-pD} \right)
\\\allowdisplaybreak&{}
=\sum_{i=0}^{\lfloor \frac{s_1}2\rfloor}\binom{s_1}{2i}\sigma_{N+i(d_1-2)}-\sum_{i=0}^{\lfloor \frac{s_1-1}2\rfloor}\binom{s_1}{2i+1} \sigma_{N+i(d_1-2)+d_1-1}\allowdisplaybreak\\&{}
=\sum_{i=0}^{\lfloor \frac{s_1}2\rfloor}\binom{s_1}{2i}\sigma_{N-2i}-\sum_{i=0}^{\lfloor \frac{s_1-1}2\rfloor}\binom{s_1}{2i+1} \sigma_{N-2i-1}\allowdisplaybreak\\&{}
=\sum _{\ell=0}^{s_1} (-1)^{\ell}\binom{s_1}{\ell}\sigma_{N-\ell},\endalign
$$for each integer $N$. Apply Lemma \tref{May9} to conclude that $\sigma_i=\sigma_j$ for all integers $i$ and $j$. 
Now that (\tref{mini}) has been established, the proof is complete.  
 \qed
\enddemo

\SectionNumber=3\tNumber=1
\heading Section \number\SectionNumber. \quad Examples.
\endheading

In Example \tref{may4} we exhibit a standard graded Artinian $k$-algebra $S$ which has an  an exact pair of homogeneous  zero divisors $(\theta_1,\theta_2)$ of degrees $2$ and $2$ without having any homogeneous exact zero divisors of degree $1$. This example is striking because the numerical result of Theorem \tref{Conj} depends on the sum $\deg \theta_1+\deg\theta_2=2+2$; but not on the particular summands $\deg \theta_1$ and $\deg\theta_2$. 

The majority of the section consists of a list of standard graded Artinian $k$-algebras which do not contain any homogeneous exact zero divisors. To obtain these rings, we apply Theorem \tref{Conj} in combination with a result from \cite{\rref{KSV}} (see Proposition \tref{KSV} below). The combined result is called Corollary \tref{comb}. Examples \tref{may5} and \tref{may6} are very explicit. In Proposition \tref{MIF*} we show that, in general, compressed level algebras do not have any homogeneous exact zero divisors. Assorted examples of compressed level algebras are given in Examples \tref{GL} (Gorenstein rings with linear resolutions), \tref{GL3} (rings defined by Pfaffians),  \tref{GL4} (determinantal rings), and \tref{Ex9} (rings corresponding to generic points in projective space).  In Proposition \tref{Mar11} and Example \tref{det}, we exhibit families of 
  standard graded Artinian $k$-algebras, which are not compressed level algebras, which nonetheless  do not contain any homogeneous exact zero divisors. These families arise from Segre embeddings and more  determinantal rings.

\example{Example \tnum{may4}}If $$S=\Bbb Q[x,y,z,w,t]/(x^4,y^4,z^4,w^4,x^2y^2,z^2w^2,y^2w^2,xt,zt,wt,t^2),$$ then $S$ has an exact pair of homogeneous zero divisors, both  of degree two, but $S$ does not have any exact zero divisors of degree one. One may check, using, for example, Macaulay2, that $\theta_1=x^2+y^2-z^2-w^2$ and $\theta_2=x^2+y^2+z^2+w^2$ is an exact pair of zero divisors. One may also check that the ideal $0:_S\ell$ is not principal for any $\ell$ of the form $\ell=a_1 x+b_2y+a_3z+a_4w+a_5t$, with $a_i\in\{0,1\}$. If $L$ is an arbitrary linear form from $S$, then there is a $\Bbb Q$-algebra automorphism of $S$ which carries $L$ to one of the $\ell$'s that have already been tested. The ideal  $0:_S\ell$ is not principal; hence the ideal  $0:_SL$ is also not principal. The Hilbert function of $S$ is
$$\matrix i&0&1&2&3&4&5&6&7&8&9\\\operatorname{HF}_S(i)&1&5&11&21&29&28&22&12&3&0.\endmatrix$$Notice that $\sigma_0=\sigma_1=\sigma_2=\sigma_3=33$, as promised by Theorem \tref{Conj} and Remark \tref{alt}. 
\endexample

\bigskip Retain the notation of Example \tref{may4} and let $$\Theta_1=x^2+y^2-z^2-w^2\quad\text{and} \quad\Theta_2=x^2+y^2+z^2+w^2$$ be pre-images, in $P=\Bbb Q[x,y,z,w,t]$, of the exact zero divisors  $\theta_1$ and $\theta_2$. Observe that the product 
$\Theta_1 \Theta_2=x^4+y^4-z^4-w^4+2x^2y^2-2z^2w^2$ is a minimal generator of the defining ideal of $S$. It is shown in \cite{\rref{KSV}} that this is a general property of exact zero divisors.

\proclaim{Proposition \tnum{KSV}} Let $P$ be a standard graded polynomial ring over the field $k$,  $I$ be a homogeneous ideal of $P$ which is primary to the maximal homogeneous ideal of $P$, and $\Theta_1$ and $\Theta_2$ be homogeneous elements of $P$ whose images in $S=P/I$ form an exact pair of homogeneous zero divisors. Then $\Theta_1\cdot \Theta_2$ is a minimal generator of $I$.\endproclaim

Throughout the rest of this section we apply Proposition \tref{KSV} in conjunction with Theorem \tref{Conj} in order to prove that various standard graded Artinian $k$-algebras do not contain any homogeneous exact zero divisors. The combined result is the following.

\proclaim{Corollary \tnum{comb}} Let $P$ be a standard graded polynomial ring over the field $k$ and  $I$ be a homogeneous ideal of $P$ which is primary to the maximal homogeneous ideal of $P$.  Assume that $[I]_1=0$,  $2\le \dim_k [P]_1$, and  $S=P/I$ contains at least  one  homogeneous exact zero divisor. 
Then the following statements hold.
\roster
\item"{(1)}" If $I$ is minimally generated by homogeneous forms of degree $n$,
 then $\frac{t^n-1}{t-1}$ divides the Hilbert series $\operatorname{HS}_S(t)$ of $S$ in the sense of {\rm(\tref{div})}. 
\item"{(2)}" If $n$ and $e$ are integers with $[I]_i=0$, for all $i$ with $i<n$, and $[I]_i=[P]_i$, for all $i$ with $e<i$, then there is an integer $D$ with $n\le D\le e$ such that  $\frac{t^D-1}{t-1}$ divides the Hilbert series $\operatorname{HS}_S(t)$ of $S$ in the sense of {\rm(\tref{div})}. 
\endroster
\endproclaim
\demo{A comment about the proof} In the situation of (2), the minimal  homogeneous  generators of $I$  have degree between $n$ and $e+1$; so a direct application of Theorem \tref{Conj}, combined with Proposition \tref{KSV}, yields that $\frac{t^D-1}{t-1}$ divides $\operatorname{HS}_S(t)$ for some $D$ with $n\le D\le e+1$. However, the ambient  hypotheses guarantee that $\dim_k[S]_0\neq \dim_k[S]_1$ and therefore $\frac{t^{e+1}-1}{t-1}$ can not possibly divide  $\operatorname{HS}_S(t)$. \qed \enddemo 

Examples \tref{may5} and \tref{may6}  involve determinantal rings. We use Corollary 1 in \cite{\rref{CH}} to determine the Hilbert functions of these rings.

\example{Example \tnum{may5}} Let $k$ be a field,  $X$ be a $4\times 5$ matrix of indeterminates, and  $I_3(X)$ be the ideal generated by the $3\times 3$ minors of $X$. Let 
$$
S=\frac{k[X]}{I_3(X)+ (\ell_1, \ldots, \ell_{14})}
$$
where $\ell_1, \ldots, \ell_{14} \in k[X]$ are linear forms such that their images in $k[X]/I_3(X)$ form a system of parameters. Then the Hilbert function of $S$ is
$$\matrix i&0&1&2&3&4&5\\\operatorname{HF}_S(i)&1&6&21&16&6&0.\endmatrix$$
If $S$ has an exact pair of homogeneous zero divisors, then the degrees of the exact zero divisors would have to add up to three according to Proposition \tref{KSV}. But note that when $D=3$, we have
$\sigma_0=17$, $\sigma_1 = 12$, $\sigma_2=21$. Therefore, Theorem \tref{Conj}, by way of Remark \tref{alt}, yields that  $S$ does not have an exact pair of homogeneous zero divisors.
\endexample

\example{Example \tnum{may6}} Let $k$ be a field,   $X$ be a $5\times 5$ matrix of indeterminates, and   $I_4(X)$ be the ideal generated by the $4\times 4$ minors of $X$. Let
$$
S=\frac{k[X]}{I_4(X)+ (\ell_1, \ldots, \ell_d)}
$$
where $\ell_1, \ldots, \ell_d\in k[X]$ are linear forms such that their images in $k[X]/I_4(X)$ form a system of parameters.
Then the Hilbert function of $S$ is
$$\matrix i&0&1&2&3&4&5&6&7\\\operatorname{HF}_S(i)&1&4&10&20&10&4&1&0.\endmatrix$$If $S$ has an exact pair of homogeneous zero divisors, then the degrees of the exact zero divisors would have to add up to four by Proposition \tref{KSV}. But note that when $D=4$, we have
$\sigma_0=11$, $\sigma_1=8$, $\sigma _2=11$, and $\sigma _3=20$. Therefore $S$ does not have an exact pair of homogeneous zero divisors. (This example is a special case of Proposition \tref{MIF*}; see Example \tref{GL4}.) In the language of Conventions \tref{D1}, $S$ an   Artinian reduction  of $\frac{k[X]}{I_4(X)}$.
\endexample

\example{Example \tnum{may7}}
Let $S$ be a standard graded Artinian algebra such that $[S]_{e+1}=0$. Let $i$ be such that $\operatorname{HF}_S(i)\ne \operatorname{HF}_S(i+1)$ and let $D > \operatorname{max}\{i+1, e-i\}$. Then $S$ cannot have any exact pair of homogeneous zero divisors with degrees adding up to $D$, since the inequality satisfied by $D$ implies that $\sigma_i=\operatorname{HF}_S(i)$ and $\sigma_{i+1}=\operatorname{HF}_S(i+1)$.
\endexample

In Propositions \tref{MIF*} and \tref{Mar11} and Example \tref{det} we give  families of standard graded Artinian $k$-algebras which do not contain any homogeneous exact zero divisors. 
The algebras of Proposition \tref{MIF*} are compressed level algebras. One introduction to this topic may be found in \cite{\rref{bmmnz}}. 
The following data is in effect.
\definition{Data \tnum{data3}} Let $S=P/I$ be a standard graded Artinian   algebra over the field $k$, where $P$ is a standard graded polynomial ring, with $c$ variables, over $k$, and $I$ is generated by homogeneous forms of degree at least two. The parameter $c$ is the {\it codimension} of $S$ in the sense that $c$ is equal to the Krull dimension of $P$ minus  the Krull dimension of $S$. \enddefinition

Retain Data \tref{data3}. Recall that the {\it socle} of $S$  is the vector space $0:_SS_+$ and the {\it type} of $S$ is $\dim_k\operatorname{socle}S$. The algebra $S$ is called {\it Gorenstein} if it has type $1$ and $S$ is called {\it level} if its socle is concentrated in  one  degree. (Notice that a Gorenstein algebra is automatically a level algebra.) Let   
$$\Bbb F:\quad 0\to \Bbb F_c@> f_c>> \Bbb F_{c-1}@>f_{c-1}>> \dots @>f_2>> \Bbb F_1@> f_1>> \Bbb F_0\tag\tnum{res}$$
 be a minimal homogeneous resolution of $S$ by free $P$-modules. If the ring $S$ is not Gorenstein, then the resolution $\Bbb F$ is called {\it linear} if the entries of $f_i$ are linear forms for $2\le i$ and the entries of $f_1$ are homogeneous forms of the same degree $n\ge 2$. If $S$ is a Gorenstein ring, then the resolution $\Bbb F$ must be symmetric; and therefore, $\Bbb F$ is called a {\it linear resolution} if  the entries of $f_i$ are linear forms for $2\le i\le c-1$ and the entries of $f_1$ and $f_c$ are homogeneous forms of the same degree $n\ge 2$.

\definition{Definition \tnum{cmped}} Retain Data \tref{data3} with $S$ a level algebra  of socle degree $e$ and type $r$. If the   Hilbert function is given by

$$\operatorname{HF}_S (i) = \min \{ \dim_k [P]_i,\ r \cdot \dim_k [P]_{e-i} \},\quad \text{for $0\le i\le e$,}$$ then $S$ is called a   {\it compressed level} algebra. \enddefinition

Notice that once the codimension, socle degree, and  type of a level algebra are fixed, then the Hilbert function of a compressed level algebra is as large as possible in each degree. 
Compressed algebras were first defined  (in a  more general context than level algebras) by   Iarrobino \cite{\rref{I84}}; where he proved that for all pairs $(e,r)$ there exists a non-empty open set of parameters which correspond to a compressed level algebra of socle degree $e$ and type $r$.     
Fr\"oberg and   Laksov \cite{\rref{FL}} offer  alternate proofs  some of Iarrobino's results. 
Zanello \cite{\rref{Za1},\rref{Za2}} has  generalized the concept to arbitrary Artinian algebras.

\proclaim{Proposition \tnum{MIF*}}  Adopt Data {\rm \tref{data3}} with $S$  a compressed level algebra with socle degree $e$, type $r$, and codimension $c$.    
If either one the following conditions hold\phantom{\tnum{mif1}\tnum{mif2}}
\roster
\item"{\rm(\tref{mif1})}" $r=1$, $2\le e$ with $e\neq 3$, and $3\le c$, or
\item"{\rm(\tref{mif2})}" $c$, $e$, and $r$ are all at least $2$ and $(c,e,r)\neq (c,2,c-1)$, 
\endroster
  then
   $S$ does not have any  homogeneous exact zero divisors. \endproclaim
\remark{\bf Remarks \tnum{e=3}} (1) The condition $e\neq 3$ is necessary   (\tref{mif1})  because $S=\frac{k[x,y,z]}{(x^2,y^2,z^2)}$ is a compressed Gorenstein algebra with  Hilbert function
$$\matrix i&0&1&2&3&4\\\operatorname{HF}_S(i)&1&3&3&1&0\endmatrix$$and socle degree $e=3$. The element $x$ of $S$ is a    homogeneous exact zero divisor.

\medskip\flushpar (2) The condition $3\le  c$ is necessary in (\tref{mif1}) because $$S=\cases \dfrac{k[x,y]}{(x^{\frac e2+1},y^{\frac e2+1})}&\text{if $e$ is even}\\
\dfrac{k[x,y]}{(x^{\frac {e+1}2},y^{\frac {e+3}2})}&\text{if $e$ is odd}
\endcases
$$ is a compressed Gorenstein algebra of socle degree $e$. The element $x$ of $S$ is a    homogeneous exact zero divisor. 

\medskip\flushpar (3) The condition  $(c,e,r)\neq (c,2,c-1)$ is necessary in (\tref{mif2}) because $$S=\frac{k[x_1,\dots, x_c]}{(x_1,\dots,x_{c-1})^2+(x_c^2)}$$ is a compressed level algebra with Hilbert function
$$\matrix i&0&1&2&3\phantom{,}\\\operatorname{HF}_S(i)&1&c&c-1&0,\endmatrix$$ codimension $c$, socle degree $e=2$, and type $c-1$. The element $x_c$ of $S$ is a    homogeneous exact zero divisor.

\endremark

 \demo{Proof of Proposition {\rm\tref{MIF*}}} Let $n$ be the minimal generator degree of $I$. We apply Corollary \tref{comb}. It suffices to show that $\frac{t^D-1}{t-1}$ does not divide $\operatorname{HS}_S(t)$, in the sense of (\tref{div}), for any integer $D$ with $n\le D\le e$. For each relevant $D$, we will exhibit two subscripts $a$ and $b$ with $\sigma_a\neq \sigma_b$ for $\sigma$ as defined in (\tref{sig}). 

Observe that $\lceil \frac {e+1} 2\rceil\le n$. Indeed, if $i<\lceil \frac {e+1} 2\rceil$, then $i\le e-i$ and 
$$\dim_k[S]_i=\operatorname{HF}_S(i)=\min\{\dim_k[P]_i,\ r\cdot \dim_k[P]_{e-i}\}=\dim_k[P]_i.$$  

If the hypotheses of (\tref{mif2}) are in effect with $e=2$, then $n=D=e=2$, the Hilbert function of $S$ is $$\matrix i&0&1&2&3\phantom{,}\\\operatorname{HF}_S(i)&1&c&r&0,\endmatrix$$ with $\sigma_0=1+r\neq c=\sigma_1$ and there is nothing more to prove. Henceforth, when the hypotheses of (\tref{mif2}) are in effect we will assume that $3\le e$.

We separate the proof into four cases. In every calculation we consider all $D$ with $$\tsize \lceil\frac {e+1}2\rceil\le n\le D\le e.$$

\medskip\flushpar{\bf Case 1.} 
Take $2\le r$, $2\le \dim_k[P]_1$,  and $4\le e$, with $e$ even. If $\frac e2+1<D$, then
$$\tsize 
\sigma_{\frac e2-1}=\operatorname{HF}_S(\frac e2-1)=\dim_k[P]_{\frac e2-1}<\dim_k[P]_{\frac e2} =\operatorname{HF}_S(\frac e2)=\sigma_{\frac e2}
,$$and if $\frac e2+1=n=D$, then
$$\split \tsize \sigma_{\frac e2-1}&\tsize =\operatorname{HF}_S(e)+\operatorname{HF}_S(\frac e2-1)=r+\dim_k[P]_{\frac e2-1}<1+r\dim_k[P]_{\frac e2-1}\\&\tsize =\operatorname{HF}_S(0)+\operatorname{HF}_S(\frac e2+1)=\sigma_{\frac e2+1}
.\endsplit $$The critical inequality holds because $2\le r$ and $2\le \dim_k[P]_{\frac e2-1}$; hence $$0<(r-1)(\dim_k[P]_{\frac e2-1}-1).$$

\medskip\flushpar{\bf Case 2.} Take $2\le r$,  $2\le \dim_k[P]_1$,  and $3\le e$, with $e$ odd. If $\frac {e+1}2<D$, then
$$\split \tsize \sigma_{\frac {e-1}2}&\tsize=\operatorname{HF}_S(\frac {e-1}2)=\dim_k[P]_{\frac {e-1}2}<\min\{\dim_k[P]_{\frac {e+1}2},\ r\cdot \dim_k[P]_{\frac {e-1}2}\}\\
&\tsize =\operatorname{HF}_S(\frac {e+1}2)=
\sigma_{\frac {e+1}2},\endsplit$$
and if $\frac {e+1}2=n=D$,
then $$\split \tsize \tsize \sigma_{\frac {e-1}2}&\tsize =\operatorname{HF}_S(\frac {e-1}2)+\operatorname{HF}_S(e)=\dim_k[P]_{\frac {e-1}2}+r<
1+r\cdot \dim_k[P]_{\frac {e-1}2}\\&\tsize =\operatorname{HF}_S(0)+\operatorname{HF}_S(\frac {e+1}2)= \sigma_{\frac {e+1}2}.
\endsplit $$

\medskip\flushpar{\bf Case 3.} Take  $r=1$,  $3\le \dim_k[P]_1$,  and  $2\le e$, with $e$  even.  If  $\frac e2+1< D\le e$, then 
$$\tsize \sigma_{\frac e2+1}=\operatorname{HF}_S(\frac e2+1)=\dim_k[P]_{\frac e2-1}<\dim_k[P]_{\frac e2}=\operatorname{HF}_S(\frac e2)=
\sigma_{\frac e2}, $$
and if $\frac e2+1= D$, then 
$$\tsize \sigma_{\frac e2+1}=\operatorname{HF}_S( 0)+\operatorname{HF}_S(\frac e2+1)=1+\dim_k[P]_{\frac e2-1}<\dim_k[P]_{\frac e2}=\operatorname{HF}_S(\frac e2)=
\sigma_{\frac e2}. $$

\medskip\flushpar{\bf Case 4.} Take  $r=1$,  $3\le \dim_k[P]_1$,  and  $5\le e$, with $e$ odd.  
Observe first that if $A$ and $B$ are integers with $1\le A$ and $2\le B$, $$\tsize\binom{A+1}1\le \binom{A+B}B\tag\tnum{duh}$$because $$\tsize A+1<\frac {A+1}1\cdot \frac {A+2}2\cdots \frac {A+B}B=\binom{A+B}B.$$ The polynomial ring $P$ has $c$ variables with $3\le c$. It follows that $$\eightpoint \dim_k [P]_1+\dim_k[P]_{B-1} <\dim_k [P]_0+\dim _k[P]_B\quad \text{for all integers $B$ with $2\le B$}\tag\tnum{use}.$$Indeed, $$\split \text{(\tref{use}) holds} &\iff \dim_k [P]_1
-\dim_k [P]_0 <\dim _k[P]_B-\dim_k[P]_{B-1}\\&\tsize \iff c-1< \binom{c-1+B}B-\binom{c-2+B}{B-1}\\&\tsize \iff \binom{(c-2)+1}1<\binom{(c-2)+B}B,\endsplit$$ and this follows from (\tref{duh}). If $\frac{e+1}2\le D\le e-2$, then apply (\tref{use}) to see that 
$$\split \sigma_1&=\operatorname{HF}_S(1)+\operatorname{HF}_S(D+1)\\&=\dim_k [P]_1+\dim_k[P]_{e-D-1}<\dim_k [P]_0+\dim _k[P]_{e-D}\\
&=\operatorname{HF}_S(0)+\operatorname{HF}_S(D)=\sigma_0.\endsplit$$
If  $D= e-1$, then 
$$\split \tsize 
\sigma_{0}&\tsize =\operatorname{HF}_S(0)+\operatorname{HF}_S(e-1)=\dim_k[P]_{0}+\dim_k[P]_{1}=c+1<\binom{c+1}2\\&\tsize=
\dim_k[P]_2=\operatorname{HF}_S(2)=\sigma_2,\endsplit $$ and if $D=e$, then 
$$\sigma_0=\operatorname{HF}_S(0)+\operatorname{HF}_S(e)
=2<c=\operatorname{HF}_S(1)=\sigma_1. \qed$$\enddemo

\proclaim{Example \tnum{GL}} Adopt Data {\rm \tref{data3}}. If $S$ is Gorenstein and the minimal homogeneous  resolution of $S$ by free $P$-modules is linear, then $S$ is a compressed level algebra. In particular, if $(c,e,r)$ satisfy {\rm(\tref{mif1})},   then $S$ does not have any homogeneous exact zero divisors. \endproclaim

\demo{Proof}The entries in the matrices $f_j$ from (\tref{res})
 are homogeneous  forms 
of degree 
$$\cases 1&\text{for $2\le j\le c-1$}\\n&\text{for $j=1$ and $j=c$},\endcases$$ for some fixed   integer $n$, with $2\le n$. In particular, $$\Bbb F_j= \cases 
P &\text{for $j=0$}\\
P(-n-j+1)^{\beta_j}&\text{for $1\le j\le c-1$}\\
P(-2n-c+2)&\text{for $j= c$}.
\endcases$$
One can use the Herzog-K\"uhl formulas \cite{\rref{HK}} to produce the betti numbers $\beta_j$, although the present argument will not use these betti numbers. One can read the socle degree $e=2n-2$ of $S$ from the back twist in $\Bbb F$; see, for example \cite{\rref{KV}, Cor.~1.7}. At this point the Hilbert function of $S$ is known. The hypothesis that the generators of $I$ have degree $n$ means that $\operatorname{HF}_S(i)=\dim_k[P]_i$, for $0\le i\le n-1$. The Hilbert function of $S$ is symmetric (because $S$ is Gorenstein) and we have already described half of the Hilbert function; therefore, the other half is known by symmetry:  
$$\operatorname{HF}_S(i)= \cases\dim_k[P]_i, &\text{for $0\le i\le n-1$}\\
\dim_k[P]_{2n-2-i}, &\text{for $n\le i\le 2n-2$}, \endcases$$ and $S$ is a compressed level algebra as described in Definition \tref{cmped}. \qed \enddemo

\proclaim{Example \tnum{GL3}} Let $n$ be an integer, $X$ be a $2n+1\times 2n+1$ alternating matrix whose entries are  linear forms from $P=k[x,y,z]$, and $I$ be the ideal in $P$ generated by the maximal order Pfaffians of $X$. If $I$ is primary to the maximal homogeneous ideal $(x,y,z)$ of $P$, then $S=P/I$ is a compressed level algebra. In particular, if   $2\le n$, then $S$ does not have any homogeneous exact zero divisors.\endproclaim
\demo{Proof} The resolution of $S$ by free $P$-modules is known \cite{\rref{BE}} to be linear. Apply Example \tref{GL}. \qed \enddemo

\remark{\bf Remark \tnum{Rem7.5}} There are plenty of compressed Gorenstein algebras which do not have linear resolutions. For example, if $S$ is described in  Data \tref{data3} with $P=k[x,y,z]$ and (\tref{res}) given by
$$0\to P(-8)\to \matrix P(-4)\\\oplus\\ P(-5)^4\endmatrix \to \matrix P(-3)^4\\\oplus\\ P(-4)\endmatrix\to P,$$
then the Hilbert function of $S$ is
$$\matrix i&0&1&2&3&4&5&6\\\operatorname{HF}_S(i)&1&3&6&6&3&1&0;\endmatrix$$hence $S$ is a compressed Gorenstein algebra which does not contain  any homogeneous exact zero divisors; but the resolution of $S$ is not linear. One such algebra $S=P/I$ is defined by $I=(x^2y, x^2z, y^3, z^3, x^4+y^2z^2)$; this ideal plays an important role in \cite{\rref{bmmnz}}.\endremark 

\bigskip The ring of Example \tref{may6} is an excellent representative of a family of compressed Gorenstein algebras which do not contain  any homogeneous exact zero divisors.
\proclaim{Example \tnum{GL4}} Let $n\ge 3$ be a positive integer, $X$ be a $(n+1)\times (n+1)$ matrix whose entries are  linear forms from $P=k[x,y,z,w]$, and $I$ be the ideal in $P$ generated by the $n\times n$ minors maximal  $X$. If $I$ is primary to the maximal homogeneous ideal $(x,y,z,w)$ of $P$, then $S=P/I$ is a compressed Gorenstein algebra, of socle degree $2n-2$, which does not contain any homogeneous exact zero divisors. \endproclaim
\demo{Proof} The resolution of $S$ by free $P$-modules is known \cite{\rref{GN}} to be linear. Apply Example \tref{GL}. \qed \enddemo

We have drawn many consequences from Proposition \tref{MIF*} when the type is one. Remark \tref{Rem8}, and Example \tref{Ex9}, are analogous to Example \tref{GL}, and Remark \tref{Rem7.5}, respectively, when the type is greater than $1$.  
\remark{\bf Remark \tnum{Rem8}} Adopt Data \tref{data3}. If $S$ is not Gorenstein and the minimal resolution of $S$ by free $P$-modules is linear, then $S$ is a compressed level algebra. Indeed, the free modules of (\tref{res}) have the form
$$\Bbb F_j=\cases P&\text{if $j=0$}\\P(-n-j+1)^{\beta_j}&\text{if $1\le j\le c$,}\endcases$$where the generators of $I$ are homogeneous forms of degree $n$. The Herzog-K\"uhl formulas \cite{\rref{HK}} give $\beta_1=\binom{n+c-1}{c-1}=\dim _k[P]_n$. It follows that $I=([P]_n)$, $S$ has socle degree $e=n-1$ and  type $r=\dim_k[P]_{n-1}$, and 
$$\operatorname{HS}_S(i)=\dim_k[P]_i=\min\{\dim_k[P]_i,\ r\cdot \dim_k[P]_{n-1-i}\},\quad \text{for $0\le i\le n-1$}.$$Thus, $S$ 
is a compressed level algebra. One could conclude that $S=\frac{k[x_1,\dots,x_c]}{(x_1,\dots,x_c)^n}$ does not contain any homogeneous exact zero divisors when $c$ and $n$ are at least two. Of course, one could also employ a direct argument to reach the same conclusion. Furthermore, $S$ is a Golod ring \cite{\rref{EG}, Prop. 1.9} of minimal multiplicity; so all of the totally reflexive $S$-modules are free; see,  \cite{\rref{AM}, (3.5)} or \cite{\rref{Y}, Cor. 2.5}. (This provides a third argument that $S$ does not contain any homogeneous exact zero divisors.)  \endremark

\example{Example \tnum{Ex9}} Let $(c,e,r)$ be a triple of positive integers which satisfy
$$\tsize \frac 1 c \binom {c+e-2}{c-1}\le r\le \binom {c+e-1}{c-1}.$$If  $n=\binom {c+e-1}{c}+r$ and $S=P/I$, as described in Data \tref{data3}, is the Artinian reduction of the coordinate ring of $n$ generic points in $\Bbb P^c$, then, according to \cite{\rref{B99}, Cor. 3.22},  $S$   is a compressed level algebra of codimension $c$, socle degree $e$, and type $r$; in particular, if $(c,e,r)$ also satisfy (\tref{mif1}) or (\tref{mif2}), then $S$ does not contain any homogeneous exact zero divisors.\endexample

\bigskip
In Proposition \tref{Mar11} and Example \tref{det} we again apply Corollary \tref{comb} to produce families of standard graded Artinian $k$-algebras which do not contain any homogeneous exact zero divisors. These families have nothing to do with compressed algebras.

 \proclaim{Proposition \tnum{Mar11}} Let $s\ge 3$ be an odd integer and let $\frak S$ be the homogeneous coordinate ring of the Segre embedding of $s$ copies of $\Bbb P^1$ into projective space. If $\ell$ is a linear system of parameters in $\frak S$ and $S=\frak S/(\ell)$, then $S$ is a standard graded Artinian Gorenstein  $k$-algebra and $S$ does not have any homogeneous exact zero divisors.\endproclaim

\demo{Proof}The Segre embedding takes $(\Bbb P^1)^s=\underbrace{\Bbb P^1\times \dots \times \Bbb P^1}_s$ into $\Bbb P^{2^s-1}$. We take $$P=k[\{W_{(j_1,\dots,j_s)}\mid j_i\in\{0,1\}\}]\quad\text {and}\quad  T=k[\{X_{i,j}\mid 1\le i\le s\text{ and }j\in\{0,1\}\}]$$ to be the coordinate rings of $\Bbb P^{2^s-1}$ and $(\Bbb P^1)^s$, respectively. The ring $\frak S$ is equal to $P/I$ for $I$ equal to the kernel of the ring homomorphism $\varphi\:P\to T$ with $\varphi(W_{j_1,\dots,j_s})=\prod\limits_{i=1}^s X_{i,j_i}$. 

The following properties of $\frak S$ are well-known; we refer to \cite{\rref{NP}}, and this reference uses results from \cite{\rref{SS}} and \cite{\rref{H}}:\phantom{\tnum{1},\tnum{2},\tnum{3}}
\roster\item"{(\tref{1})}" The ideal $I$ is generated by quadratic polynomials. 
\item"{(\tref{2})}" The ring $\frak S$ is Gorenstein of dimension $s+1$.
\item"{(\tref{3})}" The Hilbert series of $\frak S$ is $\operatorname{HS}_{\frak S}(t)=\sum\limits_{i\ge 0} (i+1)^st^i$.
\endroster
\medskip\flushpar Assertion (\tref{1}) holds because $I$ is the ``cut ideal of a tree''; see \cite{\rref{SS}, Ex. 2.3}. The cut ideal of a graph is defined in \cite{\rref{SS}}. Furthermore, in \cite{\rref{SS}}, the defining equations of the cut ideal of a graph are described in terms of the defining equations of the cut ideals of smaller graphs provided the graph $G$ may be decomposed in an appropriate manner as the ``clique sum'' of smaller graphs. Trees have the appropriate decomposition. Assertion (\tref{3}) is an immediate consequence of the fact that the Hilbert {\bf function} of a Segre product is the product of the Hilbert functions. Assertion (\tref{2}) is established in   \cite{\rref{GW}, Thms.~4.2.3, 4.4.4, 4.3.3}. In particular,   $\frak S$ is Gorenstein because it is the Segre product of rings with the same $a$-invariant. 

On the other hand, the generating function    $ \sum_{i\ge 0}(i+1)^st^i$ for $\operatorname{HS}_{\frak S}(t)$, given in (\tref{3}), is a formal power series of historical interest in Combinatorics. Indeed, if $A_s(t)$ is the polynomial defined by 
$$\sum_{i\ge 0}i^st^i=\frac{A_s(t)}{(1-t)^{s+1}},\tag\tnum{pdb}$$ then $A_s(t)$ is called the $s^{\text{th}}$ {\it Eulerian polynomial} and, if $s$ is positive, then 
$$A_s(t)=\sum_{k=1}^sA(s,k)t^k,$$ for positive integers $A(s,k)$, which are called {\it Eulerian numbers}; see  Proposition 1.4.4 and display (1.36) in \cite{\rref{St}}. The first few Eulerian polynomials are also given in \cite{\rref{St}}:
$$\matrix\format \l&\ \l &\ \l \\A_0(t)&=&1\\
A_1(t)&=&t\\
A_2(t)&=&t+t^2\\
A_3(t)&=&t+4t^2+t^3\\
A_4(t)&=&t+11t^2+11t^3+t^4\\
A_5(t)&=&t+26t^2+66t^3+26t^4+t^5.\endmatrix$$The value 
$$A_s(-1)=\cases (-1)^{(s+1)/2}E_s,&\text{if $s$ is odd}\\0,&\text{if $s$ is even and positive}\endcases$$ may be found in Exercise 135, at the end of Chapter 1 of \cite{\rref{St}}, where $E_s$  is the $s^{\text{th}}$ {\it Euler number} as defined in section 1.6.1 of \cite{\rref{St}}. The Euler numbers are positive; they satisfy the recurrence relation 
$$2E_{s+1}=\sum_{k=0}^s\binom skE_kE_{s-k}\text{ for $1\le s$, and }E_0=E_1=1;$$ and they are related to  the coefficients in the Maclaurin series for $\sec x+\tan x$ in the  sense that
$$\sec x+\tan x=\sum_{s\ge 0} E_s \frac{x^s}{s!}.$$ 

We apply Corollary \tref{comb} to prove that $S$ does not have any homogeneous exact zero divisors when   $s\ge 3 $ is an  odd integer.   The ideal $I$, which defines $S$, is generated by forms of degree $2$ according to (\tref{1}); so 
it suffices to show that $\frac{t^2-1}{t-1}=t+1$ does not  divide the Hilbert series $\operatorname{HS}_S(t)$, in the sense of (\tref{div}); that is, it suffices to show that $\operatorname{HS}_S(-1)\neq 0$. On the other hand, $\frak S$ is a Cohen-Macaulay ring of dimension $s+1$ and $S=\frak S/(\ell)$, where $\ell$ is a linear system of parameters in $\frak S$.  It follows 
that 
$$\operatorname{HS}_{\frak S}(t)=\frac{\operatorname{HS}_S(t)}{(1-t)^{\dim \frak S}}=\frac{\operatorname{HS}_S(t)}{(1-t)^{s+1}}.\tag\tnum{mr}$$
Combine (\tref{mr}), (\tref{3}), and  (\tref{pdb})   to see   that 
$$\frac{\operatorname{HS}_S(t)}{(1-t)^{s+1}}=\operatorname{HS}_{\frak S}(t)=\sum_{i\ge 0}(i+1)^st^i=\frac{\sum\limits_{i\ge 0} i ^st^i}{t}=\frac{\frac{A_s(t)}t}{(1-t)^{s+1}};$$
thus, $\operatorname{HS}_S(t)=\frac{A_s(t)}t$ and $\operatorname{HS}_S(-1)=-A_s(-1)=(-1)^{(s-1)/2} E_s\neq 0$. 
\qed
\enddemo

\proclaim{Observation \tnum{not}} If $s=3$, then the ring $S$ of Proposition {\rm\tref{Mar11}} is also studied in Proposition {\rm\tref{MIF*}}; however, if $s\ge 5$, then the ring $S$ of Proposition {\rm\tref{Mar11}} is not studied in Proposition {\rm\tref{MIF*}}. \endproclaim

\demo{Proof} If $s=3$, then the ring $S$ of Proposition \tref{Mar11} is defined by the $2\times 2$ minors of the unit cube with $W_{i_0,i_1,i_2}$ placed on the vertex $(i_0,i_1,i_2)$, with $i_j\in \{0,1\}$. The defining ideal of $S$ has nine minimal generators: one for each of the six faces of the unit cube and one for the intersection of the unit cube with each of the planes $x=y$, $x=z$, and $y=z$. The resolution of the Gorenstein ring $S$ over $P=k[\{W_{i_0,i_1,i_2}\}]$ is linear:
$$0\to P(-6)\to P(-4)^9\to P(-3)^{16}\to P(-2)^9\to P\to S,$$and $S$ is studied in Example \tref{GL}; hence also in Proposition \tref{MIF*}. (One can use Macaulay2 for this calculation.)

If $S=P/I$, as described in Data \tref{data3}, is a compressed Gorenstein algebra with the minimal generator degree of $I$ equal to $2$ and socle degree equal to $e$, then $$\min\{\dim_k[P]_2,\ \dim_k[P]_{e-2}\}=\operatorname{HF}_S(2)<\dim_k[P]_2;$$and therefore, $e\le 3$. The algebras $S$ of Theorem \tref{Mar11} are defined by ideals generated in degree $2$ and have socle degree $s-1$ because
$$\sum_{i\ge 0}\dim [S]_it^i=\operatorname{HS}_S(t)=\frac{A_s(t)}t=\sum_{i=0}^{s-1}A(s,i+1)t^i,$$ with $A(s,s)=1$. (The number $A(s,k)$ counts the permutations of $\{1,\dots,s\}$ with exactly $k-1$ descents; see \cite{\rref{St}, (1.36)}; in particular, $A(s,s)=1$.)  
If $5\le s$, then $3<s-1$; hence, the socle degree of $S$ is more than $3$ and $S$ is not a compressed Gorenstein algebra. \qed \enddemo

\proclaim{Example \tnum{det}} Let $k$ be a field, $r$ and $c$ be integers with $2\le r\le c$, $X$ be a $r\times c$ matrix of indeterminates, $I_2(X)$ be the ideal of $k[X]$ generated by the $2\times 2$ minors of $X$, and $S$ be an Artinian reduction of $k[X]/I_2(X)$, as described in Conventions {\rm (\tref{D1})}. If  any of the following conditions hold:\roster\item"{(1)}"
$r=c$ and this common number is odd, or \item"{(2)}" 
$r<c\le r+3$, or 
\item"{(3)}" $c-1$ does not divide $(r-1)!$,\endroster 
then $S$ does not have any homogeneous exact zero divisors.\endproclaim

\remark{Remark} The constraint  $c\le r+3$ in condition (2)  has been artificially imposed. We do not know any real upper constraint on $c$; indeed condition (3) applies to all large $c$, once $r$ is fixed.  \endremark

\demo{Proof} Let $R$ be the ring $k[X]/I_2(X)$ and $d$ be the Krull dimension of $R$. It is well-known; see, for example, \cite{\rref{BG}, Cor. 4},  \cite{\rref{CH}, Cor.1},  \cite{\rref{K96}}, or \cite{\rref{A88}}, that 
$$\operatorname{HS}_{R}(t)=\frac{\sum\limits_{i=0}^{r-1} \binom {r-1}i\binom {c-1}i t^i}{(1-t)^{d}}.$$ The ring $S$ is equal to $R/(\underline{\ell})$, where $\underline{\ell}$ is a regular sequence $\ell_1,\dots,\ell_d$ of homogeneous linear forms on $R$.   It follows that $\operatorname{HS}_S(t)=\sum\limits_{i=0}^{r-1} \binom {r-1}i\binom {c-1}i t^i$. We apply Corollary \tref{comb}. The defining ideal for $S$ is generated by homogeneous forms of degree $2$. It suffices to show that $\frac{t^2-1}{t-1}=t+1$ does not divide $\operatorname{HS}_S(t)$, in the sense of (\tref{div}). In other words, it suffices to show that 
$\operatorname{HS}_S(-1)\neq 0$. 

For positive integers $a\le  b$, define $N_{a,b}$ to be the integer
$$
N_{a, b}=\left\vert \sum_{i=0}^a (-1)^i \binom a   i   \binom b   i  \right\vert.
$$We see that $|\operatorname{HS}_S(-1)|=N_{r-1,c-1}$. 
Note that $N_{a, b}$ is equal (up to sign) to the coefficient of $x^a$ in $(x-1)^a(x+1)^b$, since the latter can be found as
$$
\sum _{i+j=a} (-1)^{a-i} \binom a   i  \binom b  j  
=
\sum_{j=0}^a (-1)^{j}\binom a   {a-j}   \binom b   j.  
$$
It follows from $(x-1)^a(x+1)^b =(x^2-1)^a (x+1)^{b-a}$ that
$$
N_{a, a}=\cases 0 &\text{when $a$ is odd}\\
\binom a {\frac{a}{2}}&\text{when $a$ is   even},\endcases 
$$
\medskip
$$
N_{a, a+1}=\cases\binom a { \frac{a-1}{2}}  &\text{when $a$ is odd}\\
\binom a {\frac{a}{2} }&\text{when $a$ is  even},\endcases
$$
\medskip
$$
N_{a, a+2}=\cases  2\binom a{\frac{a-1}{2}} &\text{when $a$ is  odd}\\
\binom a { \frac{a}{2}}-\binom a {\frac{a-2}{2}}&\text{when $a$ is even, and}\endcases
$$
\medskip
$$
N_{a, a+3}=\cases 3\binom a   {\frac{a-1}{2}}  -\binom a {\frac{a-3}{2}}&\text{when $a$ is  odd}\\
|\binom a {\frac{a}{2}}  -3\binom a {\frac{a-2}{2}}|   &\text{when $a$ is even}.\endcases
$$
Note that in all of the above cases (with the exception of $a=b=$ odd) we have $N_{a, b}\ne 0$.

Also note that if $N_{a, b}=0$, then $b$ must divide $a!$. In order to see this, we view $N_{a, b}$ as a polynomial of degree $a$ with rational coefficients in the variable $b$. After clearing the denominators, the constant term is $a!$ and the coefficient of $b^a$ is $(-1)^a$. \qed \enddemo


\Refs\widestnumber\key{99}


\ref\no\rnum{A88} \by S. Abhyankar \book Enumerative combinatorics of Young tableaux \bookinfo  Monographs and Textbooks in Pure and Applied Mathematics, {\bf 115} \publ Marcel Dekker, Inc. \publaddr  New York \yr  1988\endref


\ref\no \rnum{AB}  \by M. Auslander and M. Bridger 
\paper Stable module theory 
\jour Mem.   Amer. Math. Soc. \vol 94    \yr 1969 \endref  

\ref\no\rnum{AM} \by L. Avramov and A. Martsinkovsky \paper
Absolute, relative, and Tate cohomology of modules of finite Gorenstein dimension \jour
Proc. London Math. Soc. (3) \vol 85 \yr 2002   \pages 393--440\endref

\ref \no\rnum{B12} \by K. Beck \paper
Existence of totally reflexive modules via Gorenstein homomorphisms \jour
J. Commut. Algebra \vol 4 \yr 2012 \pages 57--77\endref 

\ref\no\rnum{B99} \by M.  Boij \paper Betti numbers of compressed level algebras \jour J. Pure Appl. Algebra \vol 134 \yr 1999 \pages 111--131\endref 

\ref\no\rnum{bmmnz} \paperinfo preprint, available on the arXiv\paper On the Weak Lefschetz Property for Artinian Gorenstein algebras of codimension three \by M. Boij, J. Migliore, R. M. Mir\'o-Roig, U. Nagel, and F. Zanello\endref

\ref\no\rnum{BG} \by W. Bruns and A. Guerrieri \paper The Dedekind-Mertens formula and determinantal rings \jour Proc. Amer. Math. Soc. \vol 127 \yr 1999 \pages 657--663\endref


\ref \no\rnum{BE} \by D.  Buchsbaum and D. Eisenbud \paper Algebra structures for finite free resolutions, and some structure theorems for ideals of codimension 3 \jour  Amer. J. Math. \vol 99 \yr 1977 \pages 447--485\endref

\ref\no\rnum{CJRSW} \paper Brauer-Thrall for totally reflexive modules  
\jour J.  Algebra \vol  350 \pages 340-373 \yr 2012
\by L.   W. Christensen, D. A. Jorgensen, H. Rahmati, J. Striuli, and  R. Wiegand\endref

\ref\no\rnum{CPST} \by L. W. Christensen, G. Piepmeyer, J. Striuli, and R. Takahashi \paper Finite Gorenstein representation type implies simple singularity \jour  Adv. Math. \vol 218 \yr 2008 \pages 1012--1026\endref

\ref\no\rnum{C} \by A. Conca \paper Gr\"obner bases for spaces of quadrics of low codimension \jour Adv. in Appl. Math. \vol 24 \yr 2000 \pages 111--124\endref 





\ref\no \rnum{CH} \by A. Conca and J. Herzog \paper On the Hilbert function of determinantal rings and their canonical module \jour Proc.  Amer. Math. Soc. \vol 122  \yr 1994 \pages 677--681 \endref

\ref\no\rnum{EG} \by D. Eisenbud and S. Goto, Shiro \paper
Linear free resolutions and minimal multiplicity \jour
J. Algebra \vol 88 \yr 1984 \pages 89--133\endref 


\ref\no \rnum{FL} \by R. Fr\"oberg and D. Laksov  \paper  Compressed Algebras \inbook Conference on Complete Intersections in Acireale \bookinfo Lecture Notes in Mathematics  \vol 1092 \pages 121--151 \publ Springer-Verlag \yr 1984\endref


\ref\no\rnum{GW} \by  S. Goto and K. Watanabe \paper On graded rings, I \jour J. Math. Soc. Japan \vol 30 \yr 1978) \pages  179--213\endref

\ref\no\rnum{GN} \by T. Gulliksen and O. Neg{\aa}rd \paper  Un complexe r\'esolvant pour certains id\'eaux d\'eterminantiels \jour C. R. Acad. Sci. Paris S\'er. A-B \vol 274 \yr 1972 \pages A16--A18\endref 


\ref\no \rnum{HS} \by I. Henriques and L. \c{S}ega \paper  Free resolutions over short Gorenstein local rings \jour Math. Z. \vol 267 \yr 2011\pages 645--663\endref


\ref\no \rnum{HK} \by J. Herzog, and M.  K\"uhl \paper On the Betti numbers of finite pure and linear resolutions \jour Comm. Algebra \vol 12 \yr 1984 \pages 1627--1646\endref


\ref\no\rnum{H} \by  L. T. Hoa \paper On Segre products of affine semigroup rings \jour Nagoya Math. J. \vol 110 \yr 1988 \pages 113--128\endref

\ref\no\rnum{HL} \by M. Hochster and D.  Laksov 
\paper The linear syzygies of generic forms \jour
Comm. Algebra \vol 15 \yr 1987 \pages 227--239\endref 

\ref\no\rnum{H11} \by H. Holm \paper Construction of totally reflexive modules from an exact pair of zero divisors \jour Bull. Lond. Math. Soc. \vol 43 \yr 2011 \pages 278--288\endref

\ref\no \rnum{I84} \by A. Iarrobino  \paper Compressed algebras: Artin algebras having given socle degrees and maximal length \jour Trans. Amer. Math. Soc. \vol 285  \yr 1984 \pages 337--378\endref


\ref\no \rnum{K96} \by D. Kulkarni \paper Counting of paths and coefficients of the Hilbert polynomial of a determinantal ideal \jour Discrete Math. \vol 154 \yr 1996 \pages 141--151\endref



\ref\no\rnum{KSV} \by A. Kustin, L. \c Sega, and A. Vraciu \paper Quasi-complete intersections and exact zero divisors \paperinfo in preparation \endref

\ref\no\rnum{KV} \by A. Kustin and A. Vraciu \paper Socle degrees of Frobenius powers \jour  Illinois J. Math. \vol 51 \yr 2007 \pages  185--208\endref

\ref\no \rnum{NP} \by U. Nagel and S. Petrovi\'c    \paper   Properties of cut ideals associated to ring graphs \jour J. Commut. Algebra \vol 1 \yr 2009 \pages 547--565\endref

\ref\no \rnum{SS} \by  B. Sturmfels and S. Sullivant \paper  Toric geometry of cuts and splits\jour  Michigan Math. J. \vol 57 \yr 2008 \pages 689--709\endref

\ref\no \rnum{St} \by R.  Stanley \book  Enumerative combinatorics \bookinfo Volume 1. Second edition. Cambridge Studies in Advanced Mathematics, {\bf 49} \publ Cambridge University Press \publaddr Cambridge \yr 2012\endref  


\ref\no\rnum{Y} \by Y. Yoshino \paper Modules of G-dimension zero over local rings with the cube of maximal ideal being zero \inbook Commutative algebra, singularities and computer algebra (Sinaia, 2002) \pages 255--273  \bookinfo NATO Sci. Ser. II Math. Phys. Chem. \vol 115 \publ Kluwer Acad. Publ. \publaddr Dordrecht \yr 2003\endref    

\ref\no \rnum{Za1} \by F. Zanello  \paper Extending the idea of compressed algebra to arbitrary socle-vectors \jour J. Algebra \vol 270 \yr 2003 \pages 181--198 \endref

\ref\no \rnum{Za2} \by F. Zanello  \paper Extending the idea of compressed algebra to arbitrary socle-vectors, II: cases of non-existence \jour J.  Algebra \vol 275 \yr 2004 \pages 732--750 \endref


\endRefs
\enddocument